\documentclass[11pt]{article}

\usepackage{amsmath,amssymb,enumitem,pstricks,sectsty}

\allsectionsfont{\normalsize}

\numberwithin{equation}{section}

\newtheorem{ques}{Question}

\def\BE#1{\begin{equation}\label{#1}}
\def\EE{\end{equation}}
\def\e_ref#1{(\ref{#1})}

\def\lan{\langle}
\def\ran{\rangle}
\def\lr#1{\langle{#1}\rangle}
\def\blr#1{\big\langle{#1}\big\rangle}
\def\ov#1{\overline{#1}}
\def\sf#1{\textsf{#1}}
\def\ti#1{\widetilde{#1}}

\def\wh#1{\widehat{#1}}
\def\sm#1{\begin{small}#1\end{small}}

\def\lra{\longrightarrow}

\def\De{\Delta}
\def\Ga{\Gamma}
\def\Si{\Sigma}

\def\al{\alpha}

\def\io{\iota}

\def\la{\lambda}
\def\om{\omega}

\def\i{\infty}

\def\cA{\mathcal A}
\def\C{\mathbb C}
\def\cC{\mathcal C}

\def\cH{\mathcal H}
\def\cJ{\mathcal J}
\def\fJ{\mathfrak j}
\def\L{\mathcal L}
\def\cM{\mathcal M}
\def\M{\mathfrak M}
\def\cP{\mathcal P}
\def\P{\mathbb P}
\def\Q{\mathbb Q}
\def\R{\mathbb R}
\def\V{\mathcal V}
\def\X{\mathfrak X}
\def\cZ{\mathcal Z}
\def\Z{\mathbb Z}

\def\fd{\mathfrak d}
\def\nd{\textnormal{d}}

\def\fg{\mathfrak g}
\def\fm{\mathfrak m}

\def\s{\mathbf s}

\def\nE{\textnormal{E}}

\def\nV{\textnormal{V}}

\def\Edg{\textnormal{Edg}}
\def\ev{\textnormal{ev}}
\def\Id{\textnormal{Id}}
\def\Im{\textnormal{Im}\,}
\def\PD{\textnormal{PD}}
\def\GW{\textnormal{GW}}

\def\bu{\bullet}
\def\st{\star}
\def\eset{\emptyset}

\begin{document}

\thispagestyle{empty}

\title{Some Questions in the Theory of Pseudoholomorphic Curves}
\author{Aleksey Zinger\thanks{Partially supported by NSF grant DMS 1500875}}
\date{\today}
\maketitle

\begin{abstract}
\noindent
This survey article, in honor of G.~Tian's 60th birthday, 
is inspired by R.~Pandharipande's 2002 note highlighting
research directions central to Gromov-Witten theory in algebraic geometry
and by G.~Tian's complex-geometric perspective on pseudoholomorphic curves
that lies behind many important developments
in symplectic topology since the early~1990s.
\end{abstract}

%\setcounter{section}{-1}

%\section*{Introduction}
%\label{intro_sec}

\noindent
Symplectic topology is an area of geometry originating in and closely associated 
with classical mechanics.
While long established, it has been flourishing especially  
since the introduction of pseudoholomorphic curves techniques in~\cite{Gromov}.
These techniques have led to an immense wealth of remarkable applications,
mutually enriching interplay with algebraic geometry, 
and striking connections with string theory.
They have in particular given rise to counts of such curves in symplectic manifolds,
now known as the \sf{Gromov-Witten invariants}.
While many long-standing problems have been spectacularly resolved,  
new profound questions that could have been hardly imagined in the past have arisen in their place.
This article, greatly influenced by G.~Tian's perspective on the field,
highlights a number of questions concerning pseudoholomorphic curves and their applications
in symplectic topology, algebraic geometry, and string theory.\\

\noindent
R.~Pandharipande's ICM note~\cite{Pand3ques} assembled three conjectures concerning
structures in Gromov-Witten theory:
\begin{enumerate}[label=(P\arabic*),leftmargin=*]

\item\label{cMques_it} a Poincare Duality for the tautological cohomology ring of the Deligne-Mumford
moduli space~$\ov\cM_{g,n}$ of stable nodal $n$-marked genus~$g$  curves,
known as the \sf{Gorenstein property} of $R^*(\ov\cM_{g,n})$;

\item\label{BPSques_it} integral counts of holomorphic curves in smooth complex projective threefolds,
known as the \sf{BPS~states};

\item\label{Virques_it} algebraic restrictions on Gromov-Witten invariants, known as 
the \sf{Virasoro constraints}.

\end{enumerate}
Each of these conjectures presented a deep quandary requiring fundamentally new ideas
to address.\\

\noindent
The Gorenstein property is a triviality for $g\!=\!0$, since $\ov\cM_{0,n}$ is 
a smooth projective variety and $R^*(\ov\cM_{0,n})\!=\!H^*(\ov\cM_{0,n})$.
It is established for $g\!=\!1$ in~\cite{Pe14} and shown to fail for $g\!=\!2$
whenever $n\!\ge\!20$ in~\cite{PeTo,Pe16}.
The Virasoro constraints had been established for the Gromov-Witten invariants 
of manifolds with only even-dimensional cohomology in genus~0,
of a point, of a curve, and of the complex projective space~$\P^n$
before~\cite{Pand3ques} in \cite{XiaoboTian,OP0,OP1,Giv}, respectively, 
with the last case extended to arbitrary symplectic manifolds
with semi-simple quantum cohomology in~\cite{Tel}.
However, no geometric rationale behind this conjecture that might confirm it in general
has emerged so~far, and its testing outside of fairly standard cases 
in algebraic geometry has been limited by the available computational techniques.
Just as~\ref{cMques_it},
the Virasoro Conjecture of~\ref{Virques_it} may yet turn out to fail,
at least for non-projective symplectic manifolds.
%Whenever it does hold, it may be simply a formal consequence of~\cite{Giv}
%and the pushforwards of the virtual classes of the moduli spaces of stable maps into such manifolds
%(conjecturally) being tautological classes on the moduli spaces of stable maps into~$\P^n$.
\\

\noindent
Unlike~\ref{cMques_it} and perhaps~\ref{Virques_it}, 
\ref{BPSques_it} is most naturally viewed from the symplectic topology perspective
in which it splits into three parts.
The extensive work on~\ref{BPSques_it} in algebraic geometry since~\cite{Pand3ques}
has not succeeded in confirming this conjecture even in special cases.
On the other hand, fundamentally new approaches to the three different parts of~\ref{BPSques_it}
have emerged in symplectic topology which should fully resolve 
its original formulation in a stronger formulation; see Section~\ref{BPS_sec}.\\

\noindent
The questions collected in this article fall under four distinct, but related, topics:
\begin{enumerate}[label=(\arabic*),leftmargin=*]

\item the topology of moduli spaces of pseudoholomorphic maps and applications
to the mirror symmetry predictions of string theory
and to the enumerative geometry of algebraic curves;

\item integral counts of pseudoholomorphic curves in arbitrary compact symplectic manifolds;

\item decomposition formulas for counts of pseudoholomorphic curves under
``flat" degenerations of symplectic manifolds;

\item applications of pseudoholomorphic curves techniques in symplectic topology
and algebraic geometry.

\end{enumerate}
Each of these topics involves fundamental issues concerning pseudoholomorphic curves
and a deep contribution from G.~Tian.\\

\noindent
G.~Tian's perspectives on Gromov-Witten theory had a tremendous influence
on the content of the present article in particular and the work of the author in general,
and he is very grateful to G.~Tian for generously sharing his insights
on Gromov-Witten theory over the past two decades.
The author would also like to thank J.~Li, R.~Pandharipande, and R.~Vakil
for introducing him to the richness of the algebro-geometric side of
Gromov-Witten theory indicated by many of the questions in this article
and P.~Georgieva for acquainting him with the many related mysteries of
the real sector of Gromov-Witten theory, as well as E.~Brugall\'e, 
A.~Doan, and C.~Wendl for comments on parts of this article.

\tableofcontents

\section{Topology of moduli spaces}
\label{SharpComp_sec}

\noindent
A \sf{symplectic form} on a $2n$-dimensional manifold $X$ is 
a closed 2-form on~$X$ such that $\om^n$ is a volume form on~$X$.  
A \sf{tame almost complex structure} on a symplectic manifold~$(X,\om)$ is 
a bundle endomorphism
$$J\!:TX\lra TX \qquad\hbox{s.t.}\qquad  J^2=-\Id, 
\quad \om(v,Jv)>0~~\forall\,v\!\in\!T_xX,\,x\!\in\!X,\,v\!\neq\!0.$$
If $\Si$ is a (possibly nodal) Riemann surface with complex structure~$\fJ$, 
a smooth map $u\!:\Si\!\lra\!X$ is called \sf{$J$-holomorphic} if
it solves the \sf{Cauchy-Riemann equation} corresponding to~$(J,\fJ)$:
$$\bar\partial_Ju \equiv \frac{1}{2}\big(\nd u+J\circ \nd u\circ \fJ\big) = 0.$$
The image of such a map in $X$ is called a \sf{$J$-holomorphic curve.}
\textsf{GW-invariants} are rational counts of such curves that  depend only on~$(X,\om)$.\\

\noindent
The most fundamental object in GW-theory is the moduli space $\ov\M_{g,k}(A;J)$ of 
stable $k$-marked ({\it geometric}) genus~$g$ $J$-holomorphic maps in the homology class $A\!\in\!H_2(X)$.
This compact space is generally highly singular.
However, as shown in~\cite{LT}, $\ov\M_{g,k}(A;J)$ still determines a rational homology class,
called \sf{virtual fundamental class} (\sf{VFC}) and denoted by $[\ov\M_{g,k}(A;J)]^{vir}$.
This class lives in an arbitrarily small neighborhood of $\ov\M_{g,k}(A;J)$ in
the naturally stratified configuration space $\X_{g,k}(A)$ of smooth stable maps
introduced in~\cite{LT} and is independent of~$J$.
Integration of cohomology classes against $[\ov\M_{g,k}(A;J)]^{vir}$ gives rise to GW-invariants;
see~\e_ref{GWdfn_e}.
The construction of~\cite{LT} adapts the deformation-obstruction analysis
from the algebro-geometric setting of~\cite{LT0} to symplectic topology 
via local versions of the inhomogeneous deformations the $\bar\partial_J$-equation
introduced in \cite{RT,RT2} and presents
$[\ov\M_{g,k}(A;J)]^{vir}$ as the homology class of a space stratified by 
even-codimensional orbifolds.
This approach is ideally suited for a range of concrete applications,
some of which are indicated below, and can be readily extended via~\cite{pseudo}
beyond the so-called perfect deformation-obstruction settings.\\

\begin{figure}
\begin{pspicture}(-1.1,-1.8)(10,1.25)
\psset{unit=.4cm}
\rput{45}(0,-4){\psellipse(5,-1.5)(2.5,1.5)
\psarc[linewidth=.05](5,-3.3){2}{60}{120}\psarc[linewidth=.05](5,0.3){2}{240}{300}
\pscircle[fillstyle=solid,fillcolor=gray](5,-4){1}\pscircle*(5,-3){.2}
\pscircle[fillstyle=solid,fillcolor=gray](6.83,.65){1}\pscircle*(6.44,-.28){.2}
\pscircle[fillstyle=solid,fillcolor=gray](3.17,.65){1}\pscircle*(3.56,-.28){.2}}
\rput(.2,-.9){$d_1$}\rput(3.1,2.3){$d_2$}\rput(7.8,-2.5){$d_3$}
\psarc(15,-1){3}{-60}{60}\psline(17,-1)(22,-1)\psline(16.8,-2)(21,-3)\psline(16.8,0)(21,1)
\rput(15.2,-3.5){$(1,0)$}\rput(22.4,1){$(0,d_1)$}
\rput(23.4,-1){$(0,d_2)$}\rput(22.4,-3){$(0,d_3)$}
\rput(33,-1){\begin{tabular}{l}$d_1\!+\!d_2\!+\!d_3\!=\!d$\\
$d_1,d_2,d_3\!>\!0$\end{tabular}}
\end{pspicture}
\caption{The domain of an element of $\M_1^3(\P,d)$ from the points of view 
of symplectic topology and algebraic geometry, 
with the first number in each pair in the second diagram denoting
the genus of the associated smooth irreducible component and 
the second number denoting the degree of the restriction of the map to this component.}
\label{m3_fig}
\end{figure}
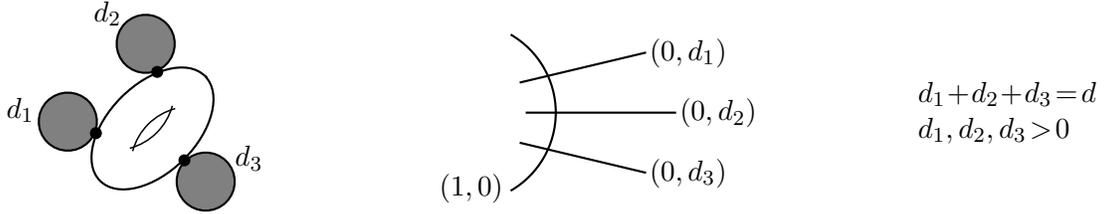

\noindent
While $\ov\M_{g,k}(A;J)$ is often called  a ``compactification"
of its subspace
$$\M_{g,k}(A;J)\subset \ov\M_{g,k}(A;J)$$
of maps from smooth domains,
$\M_{g,k}(A;J)$ usually is not a dense subset of~$\ov\M_{g,k}(A;J)$.
For example,
$$\ov\M_1(\P^n,d)\equiv \ov\M_{1,0}\big(dL;J_{\P^n}\big),$$
where $L\!\in\!H_2(\P^n)$ is the standard generator and $J_{\P^n}$ is
the standard complex structure on~$\P^n$, is a quasi-projective variety 
over~$\C$ containing $\M_1(\P^n,d)$ as a Zariski open subspace; see~\cite{FP}.
For $m\!\in\Z^+$ with $m\!\le\!n$, the dimension of the Zariski open subspace $\M_1^m(\P^n,d)$ 
of $\ov\M_1(\P^n,d)$ consisting of 
maps~$u$ from a smooth genus~1 curve~$\Si_P$ with $m$~copies of~$\P^1$ 
attached directly to~$\Si_P$ so that $u(\Si_P)\!\subset\!\P^n$ is a point~is 
$$\dim_{\C}\M_1^m(\P^n,d)=(n\!+\!1)d\!+\!n\!-\!m\ge (n\!+\!1)d=\dim_{\C}\M_1(\P^n,d);$$
see Figure~\ref{m3_fig}.
For example, 
$$\M_1^1(\P^n,d)\approx \cM_{1,1}\!\times\!\M_{0,1}(\P^n,d)\,.$$ 
Thus, $\M_1(\P^n,d)$ is not dense in $\ov\M_1(\P^n,d)$.
This motivates the following deep question concerning the convergence of
$J$-holomorphic maps in the sense of~\cite{Gromov}.

\begin{ques}[{\cite[p276]{RT}}]\label{comp_prob}
Is there a natural  Hausdorff space $\ov\M_{g,k}^0(A;J)$ of $k$~marked $J$-holomorphic maps
to~$X$ with images of arithmetic genus at least~$g$ containing $\M_{g,k}(A;J)$ as an open subspace
so that $\ov\M_{g,k}^0(A;J)$ is compact whenever $X$~is?
\end{ques}

\noindent
The ``natural" requirement in particular includes that 
$$\bigsqcup_{\begin{subarray}{c}B\in H_2(Y)\\ \io_*B=A\end{subarray}}
\!\!\!\!\!\!\ov\M_{g,k}^0\big(B;J|_Y\big)
=\big\{u\!\in\!\ov\M_{g,k}^0(A;J)\!:\Im\,u\!\subset\!Y\big\}$$
for every inclusion $\io\!:Y\!\lra\!X$ of an almost complex submanifold and
relatedly that $\ov\M_{g,k}^0(A;J)$ 
determines a fundamental class $[\ov\M_{g,k}^0(A;J)]^{vir}$.
For $g\!=\!0$, the usual moduli spaces already have the desired properties and so
$$\ov\M_{0,k}^0(A;J)=\ov\M_{0,k}(A;J)\,.$$
We also note that $\ov\M_{0,k}(\P^n,d)$ is a smooth irreducible quasi-projective variety containing 
$\M_{0,k}(\P^n,d)$ as a Zariski dense open subspace and~that 
$$\ov\M_{0,k}(\P^n,d)-\M_{0,k}(\P^n,d)\subset\ov\M_{0,k}(\P^n,d)$$
is a normal crossings divisor.\\

\noindent
For $g\!=\!1$, Question~\ref{comp_prob} is answered affirmatively in \cite{g1comp,g1comp2}
by defining
$$\ov\M_{1,k}^0(A;J)\subset\ov\M_{1,k}(A;J)$$
and showing that $\ov\M_{1,k}^0(A;J)$ determines a fundamental class.
In particular, this subspace contains an element~$u$ of $\M_{1,k}^m(A;J)$
if and only if the differentials of the restrictions of~$u$ to the $m$ copies of~$\P^1$ 
at the nodes attached to~$\Si_P$ span a subspace of $T_{u(\Si_P)}X$  
of complex dimension less than~$m$.
This imposes no condition if $2m\!>\!\dim_{\R}X$.
If $m\!\le\!n$, this imposes a condition of complex codimension $n\!+\!1\!-\!m$ 
on $\M_1^m(\P^n,d)$ and ensures that 
$$\dim_{\C}\big(\ov\M_1^0(\P^n,d)\!\cap\!\M_1^m(\P^n,d)\big)
=\dim_{\C}\M_1(\P^n,d)-1\,.$$
We also note that $\ov\M_{1,k}^0(\P^n,d)$ is a singular irreducible quasi-projective variety 
containing $\M_{1,k}(\P^n,d)$ as a Zariski dense open subspace and~that 
$$\ov\M_{1,k}^0(\P^n,d)-\M_{1,k}(\P^n,d)\subset\ov\M_{1,k}^0(\P^n,d)$$
is a divisor.
An explicit desingularization $\ti\M_{1,k}^0(\P^n,d)$ of this space is constructed 
in~\cite{g1desing} so~that 
$$\ti\M_{1,k}^0(\P^n,d)-\M_{1,k}(\P^n,d)\subset\ti\M_{1,k}^0(\P^n,d)$$
is a normal crossings divisor.
The numerical curve-counting invariants obtained by integrating cohomology 
classes against $[\ov\M_{1,k}^0(A;J)]^{vir}$ as in~\e_ref{GWdfn_e} are
called \sf{reduced genus~1 GW-invariants} in~\cite{g1comp2}.
An algebro-geometric approach to these invariants is suggested in~\cite{g1desing0}.\\

\noindent
For sufficiently positive symplectic manifolds $(X,\om)$,
the {\it standard} genus~0 and {\it reduced} genus~1 GW-invariants with insertions
pulled back from~$X$ only are {\it integer} counts of $J$-holomorphic counts
of $J$-holomorphic curves in~$X$
for a generic $\om$-compatible almost complex structure~$J$.
The standard complex structure~$J_{\P^n}$ on~$\P^n$ works for these purposes.
As demonstrated in \cite{RT,PandQdiv}, the good properties of $\ov\M_{0,k}(\P^n,d)$
indicated above are key to the enumeration of genus~0 curves in~$\P^n$
and in particular establish Kontsevich's recursion for counts of such curves.
The explicit constructions of $\ov\M_{1,k}^0(A;J)$ and $\ti\M_{1,k}^0(\P^n,d)$
in \cite{g1comp,g1desing} have opened the door for similar applications
to the enumerative geometry of genus~1 curves.\\

\noindent
For example, the Eguchi-Hori-Xiong recursion for counts of genus~1 curves in~$\P^2$
is established in~\cite{Pand99} by lifting Getzler's relation~\cite{Getz} from $\ov\cM_{1,4}$ 
to~$\ov\M_{1,k}(\P^2,d)$ and obtaining a recursion for the genus~1 GW-invariants of~$\P^2$;
the latter are the same as the corresponding enumerative invariants in this particular case.
Getzler's relation can also be lifted to $\ov\M_{1,k}(A;J)$, $\ov\M_{1,k}^0(A;J)$, and 
$\ti\M_{1,k}^0(\P^n,d)$ to yield relations between the genus~0 GW and
standard (resp.~reduced) genus~1 GW-invariants from the first (resp.~second/third) lift.
The reduced genus~1 GW-invariants of~$\P^n$ are the same as the corresponding enumerative 
invariants.
As shown in~\cite{g1diff}, the difference between the standard and reduced genus~1 GW-invariants
is a combination of the genus~0 GW-invariants;
this combination takes a very simple form in complex dimension~3.
This leads to the following, very concrete question.

\begin{ques}\label{EHX_prob} 
Can any of the above three lifts be used to obtain a recursion for the genus~1 standard or
reduced GW-invariants of~$\P^n$ for $n\!\ge\!3$ and thus a $\P^n$ analogue 
of the Eguchi-Hori-Xiong recursion  enumerating genus~1 curves?
\end{ques}

\noindent
For $g\!=\!2$, \cite{g2comp} provides the affirmative answer to
the main part of Question~\ref{comp_prob} by defining
$$\ov\M_{2,k}^0(A;J)\subset\ov\M_{2,k}(A;J)$$
and leaves no fundamental difficulty in constructing a fundamental class for this space.
The description of this subspace is significantly more complicated than
of its $g\!=\!1$ analogue.
In addition to the simple ``level~1" condition appearing in the $g\!=\!1$ case,
this description involves a more elaborate ``level~2" condition which depends
on precisely how the ``level~1" condition is satisfied relative to 
the involution and the Weierstrass points on the principal component~$\Si_P$
of the domain.
While $\ov\M_{2,k}^0(\P^n,d)$ is still a quasi-projective variety,
it is no longer irreducible and  
$\M_{2,k}(\P^n,d)$ is not dense in $\ov\M_{2,k}^0(\P^n,d)$.
However, this is not material for some applications.\\

\noindent
While Question~\ref{comp_prob} concerns a foundational issue in GW-theory
(and thus is of interest in itself),
a satisfactory answer to this problem is key to relating 
GW-invariants of a compact symplectic submanifold~$Y$ of 
a compact symplectic manifold~$(X,\om)$ given as the zero set of a transverse
bundle section  to the GW-invariants of the ambient symplectic manifold~$X$.
If $\pi_{\L}\!:\L\!\lra\!X$ is a holomorphic vector bundle and 
$\io_{\L}\!:X\!\lra\!\L$ is the inclusion as the zero section,
there is a natural projection map
\BE{Lbndl_e}
\ti\pi_{\L}\!: \V_{g,k}^A(\L)\equiv\ov\M_{g,k}(\io_{\L*}A;J)\lra\ov\M_{g,k}(A;J),
~~ \big[\ti{u}\!:\Si\!\lra\!\L\big]\lra \big[\pi_{\L}\!\circ\ti{u}\!:\Si\!\lra\!X\big].\EE
The fiber of $\ti\pi_{\L}$ over an element \hbox{$[u\!:\Si\!\lra\!X]$}
is $H^0(\Si;u^*\L)$, the space of holomorphic sections of the holomorphic 
bundle $u^*\L\!\lra\!\Si$.
If $X$ and $\L$ are sufficiently positive 
(such as $\P^n$ and sum of positive line bundles)
and $g\!=\!0$, $\ti\pi_{\L}$ is in fact a vector orbi-bundle and 
\BE{g0rel_e} \sum_{\begin{subarray}{c}B\in H_2(Y)\\ \io_*B=A\end{subarray}}
\!\!\!\!\!\!\io_*\big[\ov\M_{0,k}(B;J)\big]^{vir} 
=e\big(\V_{0,k}^A(\L)\big)\cap [\ov\M_{0,k}(A;J)]^{vir}.\EE
This observation in~\cite{Kont}, now known as the \sf{Quantum Lefschetz Hyperplane Theorem}
for genus~0 GW-invariants, was the starting point for the proofs
of the genus~0 mirror symmetry prediction of~\cite{CaDGP} 
for the quintic threefold $X_5\!\subset\!\P^4$ in~\cite{g0ms,LLY}.

\begin{ques}\label{cone_prob} 
Is there an analogue of the  $g\!=\!0$ 
Quantum Lefschetz Hyperplane Theorem~\e_ref{g0rel_e} for $g\!\ge\!1$?
\end{ques}

\noindent
While $\ti\pi_{\L}$ is not even a vector bundle for $g\!\ge\!1$
(even for sufficiently positive~$X$ and~$\L$),
it is shown in~\cite{g1cone,g1gw} that 
the restriction 
\BE{Lbndl_e2}
\ti\pi_{\L}\!:
 \V_{1,k}^A(\L)\big|_{\ov\M_{1,k}^0(A;J)} \lra  \ov\M_{1,k}^0(A;J)\EE
carries a well-defined Euler class, which in turn relates 
the reduced genus~1 GW-invariants of the submanifold~$Y$ and the ambient manifold~$X$:
\BE{g1HP_e} 
\sum_{\begin{subarray}{c}B\in H_2(Y)\\ \io_*B=A\end{subarray}}
\!\!\!\!\!\!\io_*[\ov\M_{1,k}^0(B;J)]^{vir} =
\PD_{[\ov\M_{1,k}^0(A;J)]}e\big(\V_{1,k}^A(\L)\big).\EE
This Quantum Lefschetz Hyperplane Theorem for the reduced genus~1 GW-invariants
introduced in~\cite{g1comp2} and the comparison of the standard and reduced 
genus~1 GW-invariants established in~\cite{g1diff} provide 
a Quantum Lefschetz Hyperplane Theorem for the standard genus~1 GW-invariants.
The latter is combined in~\cite{bcov1} with
the desingularization of the relevant special cases of~\e_ref{Lbndl_e2} 
constructed in~\cite{g1desing} to confirm the genus~1 mirror symmetry prediction 
of~\cite{BCOV} for~$X_5$
and to obtain similar mirror symmetry formulas for Calabi-Yau 
hypersurfaces in all projective spaces.\\

\noindent
The concrete topological construction of virtual fundamental class in~\cite{LT}
is particularly convenient for the purposes of \cite{g1comp,g1comp2,g1cone}.
It readily handles the moduli spaces $\ov\M_{1,k}^0(A;J)$,
which are not virtually smooth, but are virtually stratified by smooth orbifolds of 
even codimensions.
The representation of VFC by a geometric object in~\cite{LT} also fits well with
the comparisons carried out in \cite{g1cone,g1gw}.
However, later variations on~\cite{LT} of topological flavor, such as~\cite{FO,Pardon},
should also fit with \cite{g1comp,g1comp2,g1diff,g1cone,g1gw}.\\

\noindent
A satisfactory affirmative answer to Question~\ref{comp_prob} for each $g\!\ge\!2$,
combined with the geometric virtual fundamental class perspective of~\cite{LT},
should readily lead to a Quantum Lefschetz Hyperplane Theorem and 
to computations of GW-invariants of projective complete intersections
in the same genus~$g$.
In light of~\cite{g2comp}, there are no fundamental difficulties left to confirm
the genus~2 mirror symmetry predictions of~\cite{BCOV} for~$X_5$ and other 
projective complete intersections by paralleling the genus~1 approach initiated 
in~\cite{g1comp} and completed in~\cite{bcov1}.
The same approach should also yield confirmations of 
the mirror symmetry predictions of~\cite{Wal} for
the real GW-invariants constructed in~\cite{RealGWsI},
after the additional topological subtleties typically arising in the real setting are addressed.\\

\noindent
The methods of \cite{g1comp,g2comp} provide ``level~1" and ``level~2" obstructions
to smoothing $J$-holomorphic maps from nodal domains and can be used to define natural 
closed subspaces 
$$\ov\M_{g,k}^0(A;J)\subset\ov\M_{g,k}(A;J)$$
for $g\!\ge\!3$, which refine Gromov's Compactness Theorem and determine fundamental classes
giving rise to curve-counting invariants of compact symplectic manifolds.
However, these sharper compactifications would still not be sufficiently small 
to exclude all $J$-holomorphic maps to~$X$ with images of arithmetic genus below~$g$,
but above~1.
The associated reduced GW-invariants would then include lower-genus contributions,
even for very positive almost complex structures~$J$.
Furthermore, there are indications in~\cite{g2comp} that 
the answer to Question~\ref{comp_prob} may in fact be negative for an arbitrary almost complex
structure~$J$ on~$X$ 
if $g\!>\!2$ (or perhaps slightly larger) and the dimension of~$X$ is
sufficiently large.\\

\noindent
On the other hand, an affirmative answer to Question~\ref{comp_prob} in full generality
is not needed for specific applications, including to the enumerative geometry of positive-genus curves
in the spirit of~\cite{PandQdiv} and to the mirror symmetry predictions 
in the spirit of \cite{g1diff,bcov1}.
While the complexity of a complete description of $\ov\M_{g,k}^0(A;J)$,
whenever it can be defined,  would increase rapidly with the genus~$g$,
it is likely not to be needed for specific applications either.
In particular, it appears feasible to set up a scheme paralleling the genus~1 approach initiated 
in~\cite{g1comp} and completed in~\cite{bcov1} that would compute all GW-invariants 
of~$X_5$ modulo finitely many inputs in each genus~$g$.
This could potentially show that the generating functions~$F_g$ for these invariants
satisfy the holomorphic anomaly equations as predicted in~\cite{BCOV},
without determining each specific~$F_g$ explicitly.

\section{BPS states for arbitrary symplectic manifolds}
\label{BPS_sec}

\noindent
GW-invariants of a symplectic manifold~$(X,\om)$ are in general rational numbers
arising from families of $J$-holomorphic curves in~$X$ of possibly lower genus
and/or ``lower" degree (relative to the symplectic deformation equivalence class of~$\om$).
The \sf{primary} genus~0 GW-invariants of positive symplectic manifolds
(such as smooth Fano varieties) and of symplectic fourfolds 
arise only from $J$-holomorphic curves
of the same genus and degree, for a generic $\om$-compatible almost complex structure~$J$ on~$X$,
and are integer counts of such curves.
One might hope that the GW-invariants of~$(X,\om)$ in general are expressible
in terms of some integer invariants of~$(X,\om)$ arising from $J$-holomorphic curves
on~$X$, for $J$ generic at least in some non-empty open subset of such~$J$'s.
The explicit prediction of~\cite{GV} relating GW-invariants of 
Calabi-Yau (or~\sf{CY}) sixfolds~$(X,\om)$ to certain conjecturally integer counts
could be interpreted in such a~way;
this prediction has since been extended to a number of other special cases.\\

\noindent
For a compact symplectic manifold $(X,\om)$, 
we denote by $\cJ_{\om}$ the space of $\om$-compatible almost complex structures on~$X$.
For $g,k\!\in\!\Z^{\ge0}$, $A\!\in\!H_2(X)$, $J\!\in\!\cJ_{\om}$, and \hbox{$i\!=\!1,\ldots,k$},
let
$$\ev_i\!: \ov\M_{g,k}(A;J)\lra X, \quad \ev_i\big([u,z_1,\ldots,z_k]\big)=u(z_i),$$
be the \sf{evaluation map} at the $i$-th marked point.
We denote~by 
\BE{GWdfn_e}\begin{split}
&\qquad\GW_{g,A}^X\!: \cH^*(X)\!\equiv\!\bigsqcup_{k=1}^{\i}\!H^*(X)^{\oplus k}\lra\Q,\\
&\GW_{g,A}^X
\big(\mu_1,\ldots,\mu_k\big)= 
\bigg\lan\prod_{i=1}^k
\ev_i^*\mu_i,\big[\ov\M_{g,k}(A;J)\big]^{vir}\bigg\ran,
\end{split}\EE
the \sf{primary genus~$g$ degree~$A$ GW-invariants} of~$(X,\om)$;
these multilinear functionals are graded symmetric.
The number above vanishes unless
\BE{dimcond_e} \sum_{i=1}^k\dim_\R\mu_i
=\dim\big[\ov\M_{g,k}(A;J)\big]^{vir}=
2\big(\lr{c_1(X,\om),A}\!+\!k\big)
+\dim_{\R}X-6.\EE 
In general, this number arises from the families of genus~$g'$ degree~$A'$ $J$-holomorphic
curves in~$X$ that pass through generic pseudocycle representatives for
the Poincare duals of $\mu_1,\ldots,\mu_k$ in the sense of~\cite{pseudo}.\\

\noindent
We denote the symplectic deformation equivalence class of 
a symplectic form~$\om$ on a manifold~$X$ by~$[\om]$ and~let
$$\cA\big([\om]\big)=
\big\{(g,A)\!\in\!\Z^{\ge0}\!\times\!(H_2(X)\!-\!\{0\})\!:
\ov\M_g(A;J)\!\neq\!\eset~\forall\,J\!\in\!\cJ_{\om'},\,\om'\!\in\![\om]\big\}.$$
The genus~$g$ degree~$A$ GW-invariants of a compact symplectic manifold $(X,\om)$ 
depend only on~$[\om]$ 
and vanish unless \hbox{$(g,A)\!\in\!\cA([\om])$} or $A\!=\!0$.
In general,
they arise from families of connected $J$-holomorphic curves in~$X$ described by decorated graphs,
i.e.~tuples of the~form
\BE{Gadfn_e}\Ga=\big(\nV,\Edg,\fg\!:\nV\!\lra\!\Z^{\ge0},\fd\!:\nV\!\lra\!H_2(X)\!-\!\{0\}\big).\EE
In such a tuple, $\nV(\Ga)\!\equiv\!\nV$ and $\Edg(\Ga)\!\equiv\!\Edg$ are 
finite collections of \sf{vertices} and \sf{edges},
respectively;
the latter are pairs of vertices, but of not necessarily distinct ones, and some pairs may appear
multiple times in the collection~$\Edg$.
The vertices and the edges index the irreducible components~$\cC_v$ of the curves
and the nodes between them, respectively. 
The values of the maps $\fg$ and $\fd$  at $v\!\in\!\nV$ specify 
the geometric genus of~$\cC_v$ and its degree, respectively.
For a tuple as in~\e_ref{Gadfn_e}, we define 
$$g(\Ga)=1\!+\!|\Edg|\!-\!|\nV|\!+\!\sum_{v\in\nV}\fg(v), \quad
\fg_v(\Ga)=\fg(v),~\fd_v(\Ga)=\fd(v)~~\forall\,v\!\in\!\nV\,.$$
We denote by $\cP([\om])$ the collection of connected decorated graphs~$\Ga$ as 
in~\e_ref{Gadfn_e} such~that $(\fg(v),\fd(v))$ is an element of~$\cA([\om])$
for every $v\!\in\!\nV$.\\

\noindent
For $(g,A)\!\in\!\cA([\om])$, let $\Ga_0(g,A)$ be the unique connected edgeless graph with 
$$\fg_v\big(\Ga_0(g,A)\big)=g \qquad\hbox{and}\qquad 
\fd_v\big(\Ga_0(g,A)\big)=A$$
for the unique vertex~$v$.
Define
\begin{equation*}\begin{split}
\ti\cP_{g,A}\big([\om]\big)
=\Big\{(\Ga,\fm)\!:\Ga\!\in\!\cP([\om]),~g(\Ga)\!\le\!g,~
|\Edg(\Ga)|\!\le\!(n\!-\!3)(g\!-\!g(\Ga)),&\\
\fm\!\in\!(\Z^+)^{\nV(\Ga)},~\sum_{v\in\nV}\fm_v\fd_v(\Ga)=A&\Big\},
\end{split}\end{equation*}
where $2n\!\equiv\!\dim_{\R}X$.
By Gromov's Compactness Theorem~\cite{Gromov}, this collection 
is finite for every $(g,A)\!\in\!\cA(\om)$.
Let
$$\ti\cP_{g,A}^{\st}\big([\om]\big)\subset\ti\cP_{g,A}([\om])$$ 
be the complement of the pair $(\Ga_0(g,A),1)$.\\ 

\noindent
For a graded symmetric multilinear functional
$$\nE\!:\cH^*(X)\lra \Q$$ 
and $\mu\!\in\!H^*(X)^{\oplus k_0}$, 
we denote by $\nE(\mu,\cdot)$ the graded symmetric multilinear functional
obtained by inserting additional $k$~inputs after the $k_0$~inputs~$\mu$.
For $m\!\in\!\Z^+$, define
$$\lr{\nE}_m\!:\cH^*(X)\lra\Q, \qquad 
\lr{\nE}_m(\mu)=m^k\nE(\mu)~~\forall\,\mu\!\in\!H^k(X),\,k\!\in\!\Z^{\ge0}\,.$$ 
For graded symmetric multilinear functionals $\nE_1,\ldots,\nE_r$ as above, let
$$\prod(\nE_1,\ldots,\nE_r)\!:\cH^*(X)\lra \Q$$
be the graded symmetric multilinear functional obtained by distributing the $k$ inputs between
the $r$ functionals $\nE_1,\ldots,\nE_r$, multiplying their outputs, and 
summing over all possible distributions with the appropriate signs depending
on the degrees of the inputs.\\

\noindent
For a symplectic form $\om$ on~$X$ and $A,A^*\!\in\!H_2(X)$, we define $A\!\le_{\om}\!\!A^*$ if
$\om'(A)\!\le\!\om'(A^*)$ for some $\om'\!\in\![\om]$.
For the purposes of the question below, we identify the vertices $\nV$ of each graph
as in~\e_ref{Gadfn_e} with the set $\{1,\ldots,|\nV|\}$.

\begin{ques}\label{BPS_prob}
Let $(X,\om)$ be a compact symplectic manifold.
Are there a collection
$$C^{(g)}_{\fg,\fm}\in\Q \qquad
g\!\in\!\Z^{\ge0},~(\fg,\fm)\!\in\!(\Z^{\ge0})^r\!\times\!(\Z^+)^r,~r\!\in\!\Z^+,$$
of rational numbers and collections
$$\nE_{\Ga,\fm}^X\!:\cH^*(X)\lra\Z,~ \Ga\!\in\!\cP\big([\om]\big),\,
\fm\!\in\!(\Z^+)^{\nV(\Ga)}, \quad
\nE_{g,A}^X\!:\cH^*(X)\lra\Z,~(g,A)\!\in\!\cA([\om]),$$
of graded symmetric multilinear functionals that depend only
on~$[\om]$ and satisfy the following properties?
\begin{enumerate}[label=(E\arabic*),leftmargin=*]

\item\label{BPSform_it} for every $(g,A)\!\in\!\cA([\om])$,
\BE{BPS_e1} \GW_{g,A}^X=\nE_{g,A}^X+
\sum_{(\Ga,\fm)\in\ti\cP_{g,A}^{\st}([\om])}
\hspace{-.32in}\nE_{\Ga,\fm}^X\,;\EE

\item\label{BPSenum_it} for every $\Ga\!\in\!\cP([\om])$ as in~\e_ref{Gadfn_e},
there exist $N(\Ga)\!\in\!\Z^{\ge0}$ and $\mu_{r;v}\!\in\!\cH^*(X)$ with 
\hbox{$r\!=\!1,\ldots,N(\Ga)$} and $v\!\in\!\nV$
such~that
\BE{BPS_e1b}  \nE_{\Ga,\fm}^X=C^{(g)}_{\fg(\Ga),\fm}\!\!\sum_{r=1}^{N(\Ga)}
\prod\!\!\Big(\!\!\big(\!\blr{\nE_{\fg(v),\fd(v)}^X}_{\fm_v}\!(\mu_{r;v},\cdot)\big)_{v\in\nV}\!\Big)
\quad\forall\,\fm\!\in\!(\Z^+)^{\nV(\Ga)};\EE

\item\label{Castvan_it} for every $A\!\in\!H_2(X)$ with $\om'(A)\!>\!0$ for all $\om'\!\in\![\om]$,
\BE{BPS_e2}\sup\big\{g\!\in\!\Z^{\ge0}\!:\nE_{g,A}^X\!\neq\!0\big\}<\i;\EE

\item\label{Jenum_it} for all $g^*\!\in\!\Z^{\ge0}$ and $A^*\!\in\!H_2(X)$
there exists a subset $\cJ^{reg}_{\om}\!\subset\!\cJ_{\om}$ of second category 
in a nonempty open subset of~$\cJ_{\om}$ so that 
for all $(g,A)\!\in\!\cA([\om])$ with $g\!\le\!g^*$ and $A\!\le_{\om}\!A^*$, 
$J\!\in\!\cJ^{reg}_{\om}$, and
$\mu_1,\ldots,\mu_k\!\in\!H^*(X)$ satisfying~\e_ref{dimcond_e}, 
there exist pseudocycle representatives~$f_i$ for the Poincare duals of~$\mu_i$ 
such~that 
\begin{enumerate}[label=$\bu$,leftmargin=*]

\item the set of genus~$g$ degree~$A$ $J$-holomorphic curves meeting 
the pseudocycles $f_1,\ldots,f_k$ is cut out transversely and thus is finite, 

\item the number of such curves counted with the associated signs is 
$\nE_{g,A}^X(\mu_1,\ldots,\mu_k)$.

\end{enumerate}
\end{enumerate}
\end{ques}

\vspace{.1in}

\noindent
For all $n\!\in\!\Z^{\ge0}$ and $A\!\in\!H_2(\P^n)$,
there exists $g_A\!\!\in\!\Z^+$ so that every degree~$A$ $J_{\P^n}$-holomorphic map
\hbox{$u\!:\Si\!\lra\!\P^n$} from a smooth closed connected genus $g\!\ge\!g_A$ 
Riemann surface is a branched cover of a line $\P^1\!\subset\!\P^n$;
this is a special of the classical Castelnuovo bound \cite[p252]{GH}.
In light of~\ref{Jenum_it}, \ref{Castvan_it} is an analogue of this bound for
$J$-holomorphic curves in arbitrary symplectic manifolds.\\

\noindent
For symplectic fourfolds, i.e.~$n\!=\!2$ in the definition of the collection $\ti\cP_{g,A}([\om])$,
\e_ref{BPS_e1} and~\e_ref{BPS_e1b} reduce to $\GW_{g,A}^X\!=\!\nE_{g,A}^X$;
\ref{Jenum_it} is well-known to hold in this case.
For symplectic sixfolds, i.e.~$n\!=\!3$, $(g,A)\!\not\in\!\cA([\om])$ unless
\BE{A3cond_e} \blr{c_1(X,\om),A}=0 \qquad\hbox{or}\qquad 
\blr{c_1(X,\om),A}>0\,.\EE
In both cases, only edgeless connected graphs appear in~\e_ref{BPS_e1}.
Precise predictions for the structure of \e_ref{BPS_e1} and~\e_ref{BPS_e1b}
for symplectic sixfolds involve the coefficients $C_{h,A}(g)\!\in\!\Q$ specified~by
\BE{Ccoeff_e}\sum_{g=0}^{\i} C_{h,A}(g)t^{2g}
=\bigg(\frac{\sin(t/2)}{t/2}\bigg)^{\!2h-2+\lr{c_1(X,\om),A}}\,.\EE
In the second, \sf{Fano}, case of~\e_ref{A3cond_e}, \e_ref{BPS_e1} and~\e_ref{BPS_e1b} were 
predicted in~\cite{Pand3ques} to reduce~to 
\BE{FanoGVC_e}
\GW_{g,A}^X(\mu)=\sum_{h=0}^gC_{h,A}(g\!-\!h)\nE_{h,A}^{X}(\mu)
\qquad\forall\,\mu\!\in\!\cH^*(X).\EE
In the $g\!=\!0,1$ cases, this becomes
\BE{FanoGVC01_e}
\GW_{0,A}^X(\mu)=\nE_{0,A}^X(\mu), \quad \GW_{1,A}^X(\mu)=\nE_{1,A}^X(\mu)
+\frac{2\!-\!\lr{c_1(X,\om),A}}{24}\nE_{0,A}^X(\mu),\EE
respectively.\\

\noindent
The first equation in~\e_ref{FanoGVC01_e} with $\nE_{0,A}^X(\mu)$ described by~\ref{Jenum_it}
is the original {\it definition} of $\GW_{0,A}^X(\mu)$ for Fano classes~$A$
in the basic case of the semi-positive symplectic manifolds (which include all symplectic sixfolds).
The second equation in~\e_ref{FanoGVC01_e} holds with $\nE_{1,A}^X(\mu)$ replaced by
the reduced genus~1 GW-invariants $\GW_{1,A}^{X;0}(\mu)$ constructed in~\cite{g1comp2},
which satisfy the first bullet in~\ref{Jenum_it} whenever $(X,\om)$ is semi-positive;
see \cite[Theorem~1.1]{g1comp2} and \cite[Section~1.3]{g1comp}, respectively.
The existence of a subspace $\cJ^{reg}_{\om}\!\subset\!\cJ_{\om}$ of second category
satisfying~\ref{Jenum_it} for the Fano classes~$A$ on symplectic sixfolds is established
in~\cite{FanoGV}.
Since the system of equations~\e_ref{FanoGVC_e} with all such classes~$A$ is invertible
and the GW-invariants depend only on~$[\om]$,
this implies that the resulting counts $\nE_{h,A}^{X}(\mu)$ depend only on~$[\om]$
and thus affirmatively answers Question~\ref{BPS_prob} with the exception of~\ref{Castvan_it}
in the Fano case of~\e_ref{A3cond_e}.\\

\noindent
The first, \sf{CY}, case of~\e_ref{A3cond_e} is much harder because degree $m\!\ge\!2$ covers
of genus~$h$ degree $A/m$ $J$-holomorphic curves $\cC\!\subset\!X$ contribute to 
the genus~$g$ degree~$A$ GW-invariants of~$(X,\om)$.
For $d\!\in\!\Z^+$, we denote by $\cP(d)$ the set of partitions of~$d$
into positive integers \hbox{$d_1\!\ge\!\ldots\!\ge\!d_k$}.
Each such partition~$\rho$ corresponds to a \sf{Ferrers diagram}, 
i.e.~a collection of boxes indexed by the~set 
$$S(\rho)=\big\{(i,j)\!:i\!\in\!1,\ldots,k,~j\!\in\!1,\ldots,d_i\big\},$$
and to a \sf{dual partition} $\rho'\!\equiv\!(d_1'\!\ge\!\ldots\!\ge\!d'_{k'})$ of~$d$ specified~by
$$k'=d_1, \qquad d_j'=\max\big\{i\!=\!1,\ldots,k\!:d_i\!\ge\!j\big\}.$$
The \sf{hooklength} of a box $(i,j)\!\in\!S(\rho)$ is defined to be
$$\ell_{ij}(\rho)=d_i\!+\!d_j\!-\!i\!-\!j\!+\!1\in\Z^+\,.$$
The degree~$d$ contribution $n_{h',d}^{(h)}\!\in\!\Z^+$ of a genus~$h$ curve to the genus~$h'$ curve count
was predicted in~\cite{BryanPand} to~be given~by
\BE{CY3loc_e}\begin{split}
&\exp\Bigg(\sum_{d=1}^{\i}\sum_{h'=h}^{\i}n_{h',d}^{(h)}
\bigg(\sum_{m=1}^{\i}\frac{q^{md}}{m}\big(2\sin(m t/2)\!\big)^{\!2h'-2}\bigg)\!\!\!\Bigg)\\
&\hspace{1in}
=1+\sum_{d=1}^{\i}q^d\bigg(\sum_{\rho\in\cP(d)}\prod_{(i,j)\in S(\rho)}\!\!\!\!\!\!\!
\Big(\!2\sin\!\big(\ell_{ij}(\rho)t/2\big)\!\!\Big)^{\!2h-2}\bigg)\,.
\end{split}\EE
We note that 
\begin{gather*}
\exp\bigg(\sum_{m=1}^{\i}\frac{q^{m}}{m}\big(2\sin(m t/2)\!\big)^{-2}\bigg)
=1+\sum_{d=1}^{\i}q^d\bigg(\sum_{\rho\in\cP(d)}\prod_{(i,j)\in S(\rho)}\!\!\!\!\!\!\!
\Big(\!2\sin\!\big(\ell_{ij}(\rho)t/2\big)\!\!\Big)^{\!-2}\bigg),\\
\sum_{d=1}^{\i}\bigg(\sum_{m=1}^{\i}\frac{q^{md}}{m}\bigg)=
-\sum_{d=1}^{\i}\ln\!\big(1\!-\!q^d\big)
=\ln\!\bigg(\prod_{d=1}^{\i}\big(1\!-\!q^d\big)^{-1}\!\bigg)
=\ln\!\bigg(\!1\!+\!\sum_{d=1}^{\i}q^d\big|\cP(d)\big|\!\bigg);
\end{gather*}
the first identity above is the $t_1\!=\!t$, $t_2\!=\!t^{-1}$ case of \cite[(4.5)]{NakYos}.
Combining these two identities with the $h\!=\!0,1$ cases of~\e_ref{CY3loc_e},
we obtain
\BE{BPSloc01_e}n_{h',d}^{(0)}=\begin{cases}1,&\hbox{if}~(h',d)\!=\!(0,1);\\
0,&\hbox{otherwise};\end{cases} \quad
n_{h',d}^{(1)}=\begin{cases}1,&\hbox{if}~h'\!=\!1;\\
0,&\hbox{otherwise}.\end{cases}\EE
However, $n_{h',d}^{(h)}$ is generally nonzero for $h\!\ge\!2$, $d\!\in\!\Z^+$, 
and some $h'\!>\!h$.\\

\noindent
The primary GW-invariants~\e_ref{GWdfn_e} in the CY classes~$A$
are encoded~by the rational numbers \hbox{$N_{g,A}^X\!\equiv\!\GW_{g,A}()$},
i.e.~the GW-invariants with no insertions.
In this case, \e_ref{BPS_e1} and~\e_ref{BPS_e1b} were predicted in~\cite{GV} to
reduce~to 
\BE{CYGV_e} N_{g,A}^X=
\sum_{\begin{subarray}{c}m\in\Z^+\\ A/m\in\cA([\om])\end{subarray}}
\!\!\!\!\!\!\!\!\!m^{2g-3}
\sum_{h=0}^g\bigg(
\sum_{\begin{subarray}{c}d\in\Z^+\\ m/d\in\Z\end{subarray}}\!\!\!d^{3-2g}\!
\sum_{h'=h}^g\!C_{h',0}(g\!-\!h')n_{h',d}^{(h)}\bigg)n_{h,A/m}^{X}\,,\EE
where $n_{g,A}^X\!\equiv\!\nE_{g,A}()$.
For $m\!\in\!\Z^+$, we denote by $\lr{m}$ the sum of the positive divisors of~$m$.
By~\e_ref{BPSloc01_e}, the $g\!=\!0,1$ cases of~\e_ref{CYGV_e} become
\BE{CYGV01_e}
N_{0,A}^X=\sum_{\begin{subarray}{c}m\in\Z^+\\ A/m\in\cA([\om])\end{subarray}}
\!\!\!\!\!\!\!\!\!m^{-3}n_{0,A/m}^{X}, \quad
N_{1,A}^X=\sum_{\begin{subarray}{c}m\in\Z^+\\ A/m\in\cA([\om])\end{subarray}}
\!\!\!\!\!\!\!\!\!m^{-1}\bigg(\!\!\lr{m}n_{1,A/m}^{X}\!+\!\frac1{12}n_{0,A/m}^{X}\!\!\bigg),\EE
respectively.\\

\noindent
The system of equations~\e_ref{CYGV_e} with all CY classes~$A$ 
on a symplectic sixfold~$(X,\om)$ is also invertible.
Thus, it determines the numbers $n_{g,A}^X\!\in\!\Q$ from the number~$N_{g,A}^X$.
The original version of Question~\ref{BPS_prob}, 
known as the \sf{Gopakumar-Vafa Conjecture} for projective CY threefolds,
in fact predicted {\it only} the integrality of the numbers~$n_{g,A}^X$ obtained in this way 
and the existence of a Castelnuovo-type bound for them.
However,~\ref{Jenum_it} has been generally believed to be the underlying reason for
the validity of this conjecture since its appearance in the late 1990s. 
Until~\cite{IPgv},  \ref{Jenum_it} had also been central to every claim,
including by the authors of~\cite{IPgv} in the early~2000s, to establish the integrality
part of this conjecture;
all of these claims had quickly turned out to be erroneous.\\

\noindent
A fundamentally new perspective on the integrality part of 
the Gopakumar-Vafa Conjecture for symplectic sixfolds is introduced in~\cite{IPgv}.
It completely bypasses the analytic step~\ref{Jenum_it} and appears to succeed
in establishing the integrality of the numbers~$n_{g,A}^X$ arising from~\e_ref{CYGV_e}
via local arguments that are generally topological in spirit.
The existence of a subset $\cJ^{reg}_{\om}\!\subset\!\cJ_{\om}$ of second category
satisfying the first bullet in~\ref{Jenum_it} for symplectic CY sixfolds
is treated in~\cite{Wendl} following the general approach to this transversality
issue in~\cite{Eft}, but with additional technical input.
However, it still remains to establish that the resulting counts of $J$-holomorphic curves
satisfy the second bullet in~\ref{Jenum_it}.
Taking a geometric analysis perspective previously unexplored in GW-theory,
\cite{DoanWalpulski} uses~\cite{RiTi}, which established an analogue of Gromov's Convergence
Theorem for $J$-holomorphic maps without an a~priori genus bound,
to reduce the Castelnuovo-type bound~\ref{Castvan_it}
for symplectic CY~sixfolds to the existence of $J\!\in\!\cJ_{\om}$ 
satisfying~\ref{Jenum_it}.\\

\noindent
Precise predictions for the structure of~\e_ref{BPS_e1} and~\e_ref{BPS_e1b} 
have also been made in some cases for  symplectic manifolds of real dimensions
$2n\!\ge\!8$.
The genus~0 prediction for symplectic CY manifolds is a direct generalization
of the first equation in~\e_ref{CYGV01_e} and is given~by
\BE{CYGV01_e2}
\GW_{0,A}^X\big(\mu_1,\ldots,\mu_k\big)=
\sum_{\begin{subarray}{c}m\in\Z^+\\ A/m\in\cA([\om])\end{subarray}}
\!\!\!\!\!\!\!\!\!m^{k-3}\nE_{0,A/m}^X\big(\mu_1,\ldots,\mu_k\big)
\quad\forall\,\mu_1,\ldots,\mu_k\!\in\!H^*(X);\EE
see \cite[(2)]{KlemmP}.
The genus~1 predictions for symplectic CY manifolds of real dimensions~8 and~10 appear
in~\cite{KlemmP} and~\cite{cy5fold}, respectively.
In contrast to the arbitrary genus GW-invariants of symplectic sixfolds 
in~\e_ref{FanoGVC_e} and~\e_ref{CYGV_e} and to the genus~0 GW-invariants of symplectic CY manifolds 
in~\e_ref{CYGV01_e2}, the genus~1 GW-invariants of symplectic CY manifolds~$(X,\om)$
of real dimensions $2n\!\ge\!8$ include contributions from families of $J$-holomorphic curves
in~$(X,\om)$ of positive dimensions ($2(n\!-\!3)$-dimensional families of genus~0 curves).
This makes the analogues of~\e_ref{FanoGVC_e}, \e_ref{CYGV_e}, and~\e_ref{CYGV01_e2}  
in the last case significantly more complicated.
All curves appearing in the relevant families of $J$-holomorphic curves 
are reduced in the sense of algebraic geometric geometry and have simple nodes if $n\!=\!4,5$.
As noted in the last paragraphs of \cite[Sections~1.2,2.2]{cy5fold},
non-reduced curves and curves with non-simple nodes appear in such families if $n\!\ge\!6$.
In order to obtain a precise prediction for the structure of~\e_ref{BPS_e1} and~\e_ref{BPS_e1b}
for the genus~1 GW-invariants of symplectic CY manifolds of real dimensions $2n\!\ge\!12$,
contributions from such curves to the genus~0 and genus~1 GW-invariants
still need to be determined.\\

\noindent
Question~\ref{BPS_prob} readily extends to the real GW-invariants $\GW_{g,A}^{\phi}$
of compact real symplectic manifolds $(X,\om,\phi)$,
whenever these invariants are defined.
For example, the genus~0 real GW-invariants of real symplectic fourfolds constructed in~\cite{Wel4} 
are just signed counts of $J$-holomorphic curves.
So, the real analogues of~\e_ref{BPS_e1} and~\e_ref{BPS_e1b} in this case also
reduce to $\GW_{0,A}^{\phi}\!=\!\nE_{0,A}^{\phi}$.
Arbitrary genus real GW-invariants are constructed in~\cite{RealGWsI} for 
many real symplectic manifolds, including the odd-dimensional projective spaces~$\P^{2n-1}$
and quintic threefolds $X_5\!\subset\!\P^4$ cut out by real equations.
It is established in~\cite{NZ} that the analogue of~\e_ref{FanoGVC_e} for the Fano classes~$A$ on
a real symplectic sixfold $(X,\om,\phi)$ is 
\BE{FanoGV_e} 
\GW_{g,A}^{\phi}(\mu)=
\sum_{\begin{subarray}{c}0\le h\le g\\ g-h\in 2\Z\end{subarray}}\!\!\!\!\!
\ti{C}_{h,A}\big(\tfrac{g-h}{2}\big)\nE_{h,A}^{\phi}(\mu)
\qquad\forall~\mu\!\in\!\cH^*(X),\EE
with the coefficients $\ti{C}_{h,A}(g)\!\in\!\Q$ defined by 
\BE{Rcoeff_e}\sum_{g=0}^{\i}\ti{C}_{h,A}(g)t^{2g}
=\bigg(\frac{\sinh(t/2)}{t/2}\bigg)^{\!h-1+\lr{c_1(X,\om),A}/2}\,.\EE
The invariants $\nE_{h,A}^{\phi}(\mu)$ appearing in~\e_ref{FanoGV_e} are signed counts
of real genus~$g$ degree~$A$ $J$-holomorphic curves $\cC\!\subset\!X$.\\

\noindent
The real Fano threefold case treated in~\cite{NZ} and \cite[(5.41)]{Wal}
suggest that the real analogue of~\e_ref{CYGV_e} should~be 
\BE{CYGVR_e} N_{g,A}^{\phi}=
\sum_{\begin{subarray}{c}m\in\Z^+-2\Z\\ A/m\in\cA([\om])\end{subarray}}
\!\!\!\!\!\!\!\!\!m^{g-2}\!\!\!
\sum_{\begin{subarray}{c}0\le h\le g\\ g-h\in 2\Z\end{subarray}}\!\!\!\!\!
\bigg(\sum_{\begin{subarray}{c}d\in\Z^+\\ m/d\in\Z\end{subarray}}\!\!\!d^{2-g}
\!\!\!\!
\sum_{\begin{subarray}{c}h\le h'\le g\\ g-h'\in 2\Z\end{subarray}}\!\!\!\!\!
\ti{C}_{h',0}\big(\tfrac{g-h'}{2}\big)\ti{n}_{h',d}^{(h)}\bigg)n_{h,A/m}^{\phi}\,,\EE
for some $\ti{n}_{h',d}^{(h)}\!\in\!\Z$ (only the $d$ odd cases matter).
The right-hand sides of \cite[(5.10),(5.28)]{Wal} suggest~that 
$$\ti{n}_{h',d}^{(0)}=\begin{cases}1,&\hbox{if}~(h',d)\!=\!(0,1);\\
0,&\hbox{otherwise};\end{cases} \quad
\ti{n}_{h',d}^{(1)}=\begin{cases}1,&\hbox{if}~h'\!=\!1;\\
0,&\hbox{otherwise}.\end{cases}$$
This would reduce the $g\!=\!0,1$ cases of~\e_ref{CYGVR_e} to
$$N_{0,A}^{\phi}=\sum_{\begin{subarray}{c}m\in\Z^+-2\Z\\ A/m\in\cA([\om])\end{subarray}}
\!\!\!\!\!\!\!\!\!m^{-2}n_{0,A/m}^{\phi}, \qquad
N_{1,A}^{\phi}=\sum_{\begin{subarray}{c}m\in\Z^+-2\Z\\ A/m\in\cA([\om])\end{subarray}}
\!\!\!\!\!\!\!\!\!m^{-1}\lr{m}n_{1,A/m}^{\phi}.$$
The numbers $\ti{n}_{h',d}^{(h)}$ should arise from a real analogue of~\e_ref{CY3loc_e},
with the exponent on the left-hand side combining the real curve counts~$\ti{n}_{h',d}^{(h)}$
and the complex curve counts $n_{h',d}^{(h)}$ to account for the real doublets
of \cite[Theorem~1.3]{RealGWsII}.
The three theorems of \cite[Section~1]{RealGWsII} should provide 
the necessary geometric input to adapt the approach of~\cite{BryanPand} for~\e_ref{CYGV_e} 
to the real setting;
related equivariant localization data is provided by \cite[Section~4.2]{RealGWsIII}.
Some analogue of~\e_ref{CY3loc_e} for the real setting has been apparently 
obtained by~\cite{RealGV}.\\

\noindent
The approach of \cite{IPgv} to the integrality of the numbers~$n_{g,A}^X$ determined
by~\e_ref{CYGV_e} should be adaptable to other situations when the GW-invariants in question
 are expected
to arise entirely from isolated $J$-holomorphic curves.
These situations include the real genus~0 GW-invariants of many real symplectic manifolds 
and the real arbitrary genus GW-invariants of real symplectic CY sixfolds constructed
in~\cite{Ge2} and~\cite{RealGWsI}, respectively.
In fact, the integrality of the numbers~$\nE_{0,A}^X(\mu)$ determined
by~\e_ref{CYGV01_e2} is already a (secondary) subject of~\cite{IPgv}.
On the other hand, the approach of~\cite{IPgv} does not appear readily adaptable
to situations when positive-dimensional families of $J$-holomorphic curves in~$X$
are expected to contribute to the GW-invariants in question.
These situations include the genus~1 GW-invariants of symplectic CY manifolds
of real dimensions~8 and~10 studied in~\cite{KlemmP} and~\cite{cy5fold}, respectively.
The approaches of~\cite{Wendl} and~\cite{DoanWalpulski} 
to the existence of a subset $\cJ^{reg}_{\om}\!\subset\!\cJ_{\om}$
satisfying~\ref{Jenum_it} and to the Castelnuovo-type bound for the associated counts
of $J$-holomorphic curves, respectively, appear more flexible in this regard.\\

\noindent
Enumerative geometry of curves in projective varieties is a classical subject 
originating in the middle of the nineteenth century.
However, the developments in this field had been limited to very low degrees until
the emergence of GW-theory and its applications to enumerative geometry in the early~1990s.
As the moduli spaces $\ov\M_{g,k}(A;J)$ have fairly nice deformation-obstruction theory,
the GW-invariants arising from these spaces are often amendable to computations.
Whenever these invariants can be related to enumerative curve counts as in Question~\ref{BPS_prob},
computations of GW-invariants translate into direct applications to enumerative geometry.
The most famous such application is perhaps Kontsevich's recursion for counts
of genus~0 curves in~$\C\P^2$, stated in~\cite{KM} and proved in~\cite{RT}.
Analogues of this recursion for counts of real genus~0 curves in~$\P^2$ defined in~\cite{Wel4}
and in~$\P^{2n-1}$ defined in~\cite{Ge2} appear in~\cite{Sol2} and
\cite{RealEnum,RealEnumApp}, respectively.
The counts of genus~$g$ degree~$d$ curves arising from
the proofs of the mirror symmetry predictions for the projective CY complete intersections 
in genus~0 in~\cite{g0ms,LLY} and in genus~1 in~\cite{bcov1,Po} via~\e_ref{CYGV01_e} 
have been shown to match the classical enumerative counts for $g\!=\!0$, $d\!\le\!3$
and for $g\!=\!1$, $d\!\le\!4$; see~\cite{ES}.
The genus~0  real GW-invariants of real symplectic fourfolds defined in~\cite{Wel4}
and of many higher-dimensional real symplectic manifolds defined in~\cite{Ge2}
directly provide lower bounds for counts of genus~0 real curves;
the arbitrary genus real GW-invariants defined in~\cite{RealGWsI}
provide such bounds in arbitrary genera via the relation~\e_ref{FanoGV_e} 
proved in~\cite{NZ}.
For local CY manifolds, Question~\ref{BPS_prob} points to intriguing number-theoretic properties
of GW-invariants; see G.~Martin's conjecture in \cite[Section~3.2]{cy5fold}.\\

\noindent
The coefficients $\ti{C}_{0,0}(g)$ in~\e_ref{Rcoeff_e} are the coefficients 
of the renown $A$-series central to the Index Theorem \cite[Theorem~3.13]{LM};
they in particular determine the index of the Dirac operator on a Spin~bundle.
The coefficients in~\e_ref{Ccoeff_e} are closely related to the $A$-series as well.
It is tempting to wonder if there is some connection between the multiply covered contributions
encoded by~\e_ref{Ccoeff_e} and by~\e_ref{Rcoeff_e} and Dirac operators.

\section{Symplectic degenerations and Gromov-Witten invariants}
\label{GWdecomp_sec}

\begin{figure}
\begin{pspicture}(-1.5,-1)(11,2)
\psset{unit=.3cm}
\psline[linewidth=.1](5,-2)(12,-2)\psline[linewidth=.1](5,-2)(5,5)
\pscircle*(5,-2){.3}\psarc(11,4){5}{180}{270}
\rput(10,1.5){\sm{$X_1\!\#_{X_{12}}\!X_2$}}
\rput(11.5,-2.9){\sm{$X_1$}}\rput(4,4.5){\sm{$X_2$}}
\rput(4.5,-2.8){\sm{$X_{12}$}}
%2nd diagram
\psline[linewidth=.1](20,0)(28,0)
\pscircle*(20,0){.3}\pscircle*(28,0){.3}
\rput(19.2,-.8){\sm{$X_1$}}\rput(28.8,-.8){\sm{$X_2$}}
\rput(24,-.8){\sm{$X_{12}$}}
% 3rd diagram
\psline[linewidth=.1](35,-3)(43,-3)\psline[linewidth=.1](35,-3)(35,5)
\psline[linewidth=.1](43,-3)(35,5)
\psline[linewidth=.1](35,0)(40,0)\pscircle*(40,0){.3}
\pscircle*(35,-3){.3}\pscircle*(43,-3){.3}\pscircle*(35,5){.3}
\pscircle*(35,0){.3}
\rput(37,-.8){\sm{$E$}}\rput(37,.8){\sm{$L$}}
\end{pspicture}
\caption{A 2-fold simple normal crossings variety $\cZ_0\!=\!X_1\!\cup_{X_{12}}\!X_2$
with its smoothing \hbox{$\cZ_{\la}\!=\!X_1\!\#_{X_{12}}\!X_2$}, 
its dual intersection complex,
and a toric 2-fold decomposition of~$\P^2$ into $\P^2$ and its one-point blowup $\wh\P^2$
along a line \hbox{$L\!\subset\!\P^2$} and the exceptional divisor \hbox{$E\!\subset\!\wh\P^2$}.}
\label{SympSum2_fig}
\end{figure}
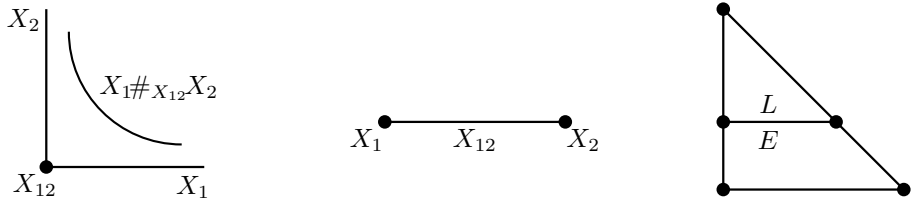

\noindent
It is natural and essential to study the behavior of GW-invariants under reasonable
degenerations and decompositions of symplectic manifolds, as pointed out in~\cite{Ti}.
The standard example of such a decomposition is provided by 
the \sf{symplectic sum construction} of~\cite{Gf};
it joins two symplectic manifolds~$X_1$ and~$X_2$ along a common {\it smooth} 
\sf{symplectic divisor}~$X_{12}$
(i.e.~a closed symplectic submanifold of real codimension~2)
with dual normal bundles in the two manifolds into a symplectic manifold~$X_1\!\#_{X_{12}}\!X_2$.
In fact, the symplectic sum construction of~\cite{Gf} provides 
a symplectic fibration \hbox{$\pi\!:\cZ\!\lra\!\De$} over the unit disk $\De\!\subset\!\C$, 
whose central fiber~$\cZ_0$ is $X_1\!\cup_{X_{12}}\!X_2$ and the remaining fibers are 
smooth symplectic manifolds which are symplectically deformation equivalent to each other;
see the first diagram in Figure~\ref{SympSum2_fig}.
While the behavior of GW-invariants under the basic degenerations and decompositions
associated with the construction of~\cite{Gf} was understood long ago and has since been
followed up by numerous applications throughout GW-theory,
the progress beyond these cases has been slow.
The interest in finding usable decomposition formulas for GW-invariants in more general situations
has grown considerably since the advent of the \sf{Gross-Siebert program}~\cite{GS0}
for a (fairly) direct approach to the mirror symmetry predictions of string theory.\\

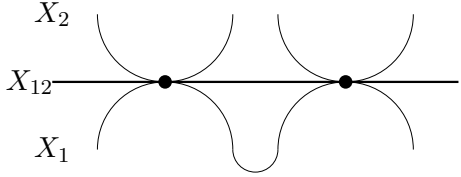
\begin{figure}
\begin{pspicture}(38,-2.3)(11,1)
\psset{unit=.3cm}
% 2nd diagram starts here
\psline[linewidth=.1](55,-2)(73,-2)
\rput(54,-2){$X_{12}$}\rput(55,-5){$X_1$}\rput(55,1){$X_2$}
\psarc[linewidth=.04](60,-5){3}{0}{180}\pscircle*(60,-2){.3}
\psarc[linewidth=.04](64,-5){1}{180}{0}
\psarc[linewidth=.04](68,-5){3}{0}{180}\pscircle*(68,-2){.3}
\psarc[linewidth=.04](60,1){3}{180}{0}\psarc[linewidth=.04](68,1){3}{180}{0}
\end{pspicture}
\caption{A connected curve in~$\cZ_0$ possibly contributing to the GW-invariants of~$\cZ_{\la}$.}
\label{limitcurve_fig}
\end{figure}

\noindent
A sequence of $J$-holomorphic curves in the smooth fibers $\cZ_{\la}\!=\!X_1\!\#_{X_{12}}\!X_2$ 
of a symplectic fibration \hbox{$\pi\!:\cZ\!\lra\!\De$} associated with the construction
of~\cite{Gf} with $\la\!\lra\!0$ converges to curves in the singular fiber 
$\cZ_0\!=\!X_1\!\cup_{X_{12}}\!X_2$.
Each of the irreducible components of a limiting curve 
either lies entirely in~$X_{12}$ or meets~$X_{12}$ in finitely many points (possibly none)
and lies entirely in either $X_1$ or~$X_2$.
A key prediction in~\cite{Ti} concerning the behavior of the GW-invariants of~$\cZ_{\la}$
as $\la\!\lra\!0$ is that they should arise only from $J$-holomorphic curves in~$\cZ_0$ 
with no irreducible components contained in~$X_{12}$ and with the irreducible components mapped
into~$X_1$ and~$X_2$ having the same contacts with~$X_{12}$; see Figure~\ref{limitcurve_fig}.
In particular, there should be {\it no} direct contribution from the GW-invariants of~$X_{12}$.
The multiplicity with which such a limiting curve should contribute to
the GW-invariants of~$\cZ_{\la}$ is determined in~\cite{CH} based on 
a straightforward algebraic reason.\\

\noindent
Notions of stable $J$-holomorphic maps to  
\sf{simple normal crossings} (or~\sf{SC}) projective varieties of
the form $X_1\!\cup_{X_{12}}\!X_2$ and of stable maps to~$X_i$ relative to
a smooth projective divisor~$X_{12}$ are introduced in~\cite{JunLi1}.
A \sf{degeneration formula} relating 
the virtual cycles of the moduli spaces $\ov\M_{g,k}(A_{\la};J)$ with 
$A_{\la}\!\in\!H_2(\cZ_{\la})$ to 
the virtual cycles  of the moduli spaces $\ov\M_{g,k}(A_0;J)$ with $A_0\!\in\!H_2(\cZ_0)$ 
appears in~\cite{JunLi2}.
A \sf{splitting formula} decomposing 
the latter  into the virtual cycles of the moduli spaces 
$\ov\M_{g_i,k_i;\s_i}(X_{12},A_i;J)$ of stable relative maps to~$(X_i,X_{12})$ via 
a Kunneth decomposition of the~diagonal
$$\De_{X_{12}}=\big\{(x,x)\!:x\!\in\!X_{12}\big\}$$
in $X_{12}^{\,2}$ is also established in~\cite{JunLi2}.
The relative GW-invariants of~$(X_i,X_{12})$ are in turn shown to 
reduce to the (absolute) GW-invariants of~$X_i$ and~$X_{12}$ in~\cite{MauP}.
Thus, \cite{JunLi1,JunLi2,MauP} fully establish the prediction of~\cite{Ti}
in the projective category in the case of basic degenerations of the target
as in Figure~\ref{SympSum2_fig}.
An expository account of the symplectic topology perspective 
on the numerical reduction of the decomposition formula of~\cite{JunLi2}
appears in~\cite{SympSum}.\\

\noindent
The standard symplectic sum construction of~\cite{Gf} readily extends to the setting
where the disjoint union $X_1\!\sqcup\!X_2$ is replaced by a single symplectic manifold
$(\ti{X},\ti\om)$ and 
the two copies of the divisor~$X_{12}$ are replaced by 
a single smooth symplectic divisor $\ti{X}_{12}\!\subset\!\ti{X}$  
with a symplectic involution~$\psi$.
The \sf{NC symplectic variety} \hbox{$\cZ_{\psi;0}\!\equiv\!X_{\psi}$} is 
then obtained from~$\ti{X}$ by identifying
the points on~$\ti{X}_{12}$ via~$\psi$.
This setting is discussed in Example~6.10 in the first two versions of~\cite{SympDivConf};
a construction smoothing $\cZ_{\psi;0}$ 
into symplectic manifolds $\cZ_{\psi;\la}$
is a special case of the construction outlined in Section~7 of 
the first version of~\cite{SympSumMulti} and detailed in~\cite{SympSumMulti2}.
The reasoning behind the decomposition formulas for GW-invariants in
the basic setting of the previous paragraph readily extends to provide a relation
between the GW-invariants of a smoothing~$\cZ_{\psi;\la}$ of the NC symplectic variety~$X_{\psi}$ 
and the relative GW-invariants of~$(\ti{X},\ti{X}_{12})$.
The only difference in the resulting formula is that a Kunneth decomposition
of the~diagonal \hbox{$\De_{X_{12}}\!\subset\!X_{12}^{\,2}$}
is replaced by a Kunneth decomposition of the $\psi$-diagonal
$$\ti\De_{\psi}=\big\{\big(\ti{x},\psi(\ti{x})\!\big)\!:\ti{x}\!\in\!\ti{X}_{12}\big\};$$
the resulting sum of pairwise products of the GW-invariants of $(\ti{X},\ti{X}_{12})$
should then be divided by~2.\\

\noindent
The decomposition formulas of \cite{JunLi2} do not completely determine 
the GW-invariants of a smooth fiber $\cZ_{\la}\!=\!X_1\!\#_{X_{12}}\!X_2$
in terms of the GW-invariants of $(X_i,X_{12})$ in many cases because of 
the so-called \sf{vanishing cycles}:
second homology classes in~$\cZ_{\la}$ which vanish under the projection 
to \hbox{$\cZ_0\!=\!X_1\!\cup_{X_{12}}\!X_2$}.
A refinement to the usual relative GW-invariants of~$(X,V)$ of~\cite{JunLi1} 
is suggested in~\cite{IPrel} with the aim of resolving this unfortunate deficiency of 
the decompositions formulas of~\cite{JunLi2} in~\cite{IPsum}.
This refinement is constructed in~\cite{GWrelIP} via a lifting 
$$\ti\ev_X^V\!:\ov\M_{g,k;\s}(V,A;J)\lra \wh{V}_{X;\s}$$
of the relative evaluation map to a covering of $V_{\s}\!\equiv\!V^{\ell}$,
where $\ell\!\in\!\Z^{\ge0}$ is the length of the relative contact vector~$\s$.
This refinement sharpens the decomposition 
formulas of~\cite{JunLi2} by pulling back closed submanifolds
\BE{diagsplit_e}\wh{V}_{X_1,X_2;\s}^A\subset 
\big(\wh{V}_{X_1;\s}\!\times\!\wh{V}_{X_2;\s}\big)\big|_{\De_V^{\,\ell}}\,,\EE
with $V\!=\!X_{12}$ and $A\!\in\!H_2(X_1\!\#_{X_{12}}\!X_2)$, by 
$\ti\ev_{X_1}^V\!\times\!\ti\ev_{X_2}^V$; see \cite[Section~1.2]{GWsumIP}.
However, this does not necessarily lead to a decomposition of 
the GW-invariants of~$X_1\!\cup_{X_{12}}\!X_2$ into the GW-invariants 
of~$(X_i,X_{12})$ that completely describes the former in terms of the latter.
The same approach provides a sharper version of the relation between
the GW-invariants of a smoothing~$\cZ_{\psi;\la}$ of~$X_{\psi}$ and 
the GW-invariants of $(\ti{X},\ti{X}_{12})$ indicated in the previous paragraph.
The submanifolds~\e_ref{diagsplit_e} in this case are replaced by certain submanifolds
$$\wh{V}_{\ti{X};\s\s}^A\subset \wh{V}_{\ti{X};\s\s}\big|_{\De_{\psi}^{\,\ell}}\,,$$
with $V\!=\!\ti{X}_{12}$;
the resulting  relative invariants of $(\ti{X},\ti{X}_{12})$ should then be divided by~2.\\

\noindent
{\it Qualitative} applications of the above refinements to relative GW-invariants and 
to the decomposition formula of~\cite{JunLi2} are described in \cite{GWrelIP,GWsumIP}.
These refinements in principle distinguish between the GW-invariants of~$\cZ_{\la}$
in degrees~$A_{\la}$ differing by torsion.
Torsion classes can also arise from the one-parameter families of smoothings~$\cZ_{\psi;\la}$
of~$X_{\psi}$ as above.
{\it Quantitative} computation of GW-invariants in degrees differing by torsion has been 
a long-standing problem.

\begin{ques}\label{torsionGW_prob}
Is it possible to compute GW-invariants in degrees differing by torsion in some cases
via the sharper version of the decomposition formula described in~\cite{GWsumIP}
and/or its analogue for the degenerations of the form~$\cZ_{\psi;\la}$ above?
\end{ques}

\noindent
The Enriques surface~$X$ forms an elliptic fibration over~$\P^1$ with 
12 nodal fibers and 2~double fibers; see \cite[Section~1.3]{MauP2}.
The difference $F_1\!-\!F_2$ between the two double fibers is a 2-torsion class.
A~smooth genus~1 curve~$E$ has a fixed-point-free holomorphic involution~$\psi$.
The quotient
$$X_2\equiv \big(\P^1\!\times\!E\big)\big/\!\!\sim, \qquad
(z,p)\lra \big(\!-\!z,\psi(p)\big),$$
forms an elliptic fibration over~$\P^1$ with 2~double fibers.
The blowup $\ti{X}$ of $\P^2$ at the 9-point base locus of a generic pencil of
cubics is an elliptic fibration over~$\P^1$ with 12~nodal fibers.
The NC variety $\cZ_0\!\equiv\!X_2\!\cup_E\!\ti{X}$ can be smoothed out to 
an Enriques surface $\cZ_{\la}\!\equiv\!X$.
The genus~1 GW-invariants of~$X$ are determined in~\cite{MauP2}
by applying the decomposition formula of~\cite{JunLi2} 
in this setting and using the Virasoro constraints.
However, the computation in~\cite{MauP2} does not distinguish between the map degrees differing
by the torsion $F_1\!-\!F_2$;
this torsion arises from the vanishing cycles and thus is not detected by 
the decomposition formula of~\cite{JunLi2}.
On the other hand, it may be possible to fully compute the genus~1 GW-invariants of~$X$ 
by refining the computation in~\cite{MauP2} via
the sharper version of this formula described in~\cite{GWsumIP}.\\

\noindent
Another potential approach to a complete computation of the GW-invariants of
the Enriques surface~$X$ is provided by the extension of the standard symplectic sum
construction of~\cite{Gf} indicated above Question~\ref{torsionGW_prob}.
Let $\psi$ be a fixed-point-free holomorphic involution on a smooth fiber $F\!\approx\!E$
of \hbox{$\ti{X}\!\lra\!\P^1$}.
The NC variety
$$\cZ_{\psi;0}\!\equiv\!X_{\psi}\equiv \ti{X}\!\big/\!\!\sim, \qquad
p\sim\psi(p)~~\forall~p\!\in\!\ti{X}_{12}\!\equiv\!F,$$
has a $\Z^2$-collection of one-parameter families of smoothings~$\cZ_{\psi;\la}$.
The total spaces of these families are $\Z_2$-quotients of 
the total families of the smoothings of $\ti{X}\!\cup_F\!\ti{X}$.
The fibers~$\cZ_{\la}$ of one of the latter families are $K3$ surfaces.
Thus, the fibers~$\cZ_{\psi;\la}$ in one of the families of smoothings of~$X_{\psi}$ should
be Enriques surfaces (at least up to symplectic deformation equivalence). 
The extension of the standard degeneration formula of~\cite{JunLi2} indicated above
applies to these families of smoothings and again distinguishes between the GW-invariants
in degrees differing by the torsion $F_1\!-\!F_2$.\\

\begin{figure}
\begin{pspicture}(-2,-2)(11,2)
\psset{unit=.3cm} 
\psline[linewidth=.1](5,-2)(12,-2)\psline[linewidth=.1](5,-2)(5,5)
\psline[linewidth=.1](5,-2)(0.5,-6.5)\pscircle*(5,-2){.3}
\rput(9.5,2.5){\sm{$X_1$}}\rput(1.5,-1){\sm{$X_2$}}\rput(7,-5){\sm{$X_3$}}
\rput(11.5,-1.1){\sm{$X_{13}$}}\rput(6.2,4.5){\sm{$X_{12}$}}\rput(-.4,-5.7){\sm{$X_{23}$}}
\rput(6.5,-1.2){\sm{$X_{123}$}}
%2nd diagram
\psline[linewidth=.1](20,-4)(20,4)\psline[linewidth=.1](20,-4)(28,-4)
\psline[linewidth=.1](28,-4)(20,4)
\pscircle*(20,-4){.3}\pscircle*(20,4){.3}\pscircle*(28,-4){.3}
\rput(22.5,-1.5){\sm{$X_{123}$}}
\rput(19.2,-4.8){\sm{$X_1$}}\rput(28.8,-4.8){\sm{$X_2$}}\rput(19.2,4.8){\sm{$X_3$}}
\rput(24,-4.8){\sm{$X_{12}$}}\rput(18.8,0){\sm{$X_{13}$}}\rput(25.5,0){\sm{$X_{23}$}}
% 3rd diagram
\psline[linewidth=.1](35,-5)(45,-5)\psline[linewidth=.1](35,-5)(35,5)
\psline[linewidth=.1](45,-5)(35,5)
\psline[linewidth=.1](35,-1)(38,-1)\psline[linewidth=.1](38,2)(38,-1)
\psline[linewidth=.1](42,-5)(38,-1)\pscircle*(38,-1){.3}
\pscircle*(35,-5){.3}\pscircle*(45,-5){.3}\pscircle*(35,5){.3}
\pscircle*(35,-1){.3}\pscircle*(42,-5){.3}\pscircle*(38,2){.3}
\rput(40.6,-2.4){\sm{$E$}}\rput(39.3,-3.4){\sm{$\bar{L}$}}
\rput(37.2,.3){\sm{$E$}}\rput(38.5,.3){\sm{$\bar{L}$}}
\rput(35.8,-1.8){\sm{$E$}}\rput(35.8,-.2){\sm{$\bar{L}$}}
\end{pspicture}
\caption{A 3-fold simple normal crossings variety~$\cZ_0$, 
its dual intersection complex,
and a toric 3-fold decomposition of~$\P^2$ into three copies of its one-point blowup $\wh\P^2$
 along the exceptional divisor $E\!\subset\!\wh\P^2$ and the proper transform
$\ov{L}\!\subset\!\wh\P^2$ of a line $L\!\subset\!\P^2$.}
\label{SympSum3_fig}
\end{figure}
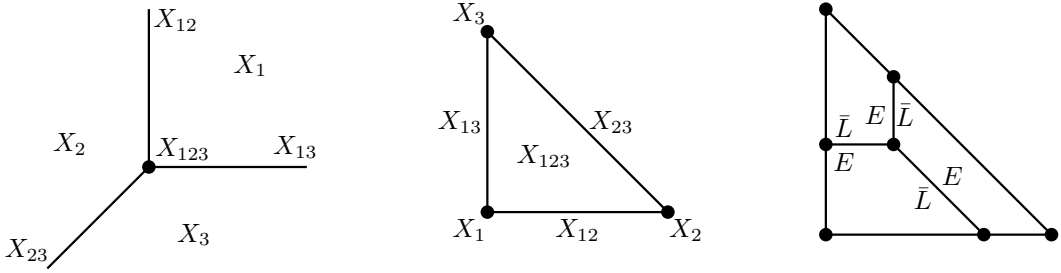

\noindent
The Gross-Siebert program~\cite{GS0} for a direct proof of mirror symmetry 
requires degeneration and splitting formulas for GW-invariants under degenerations
\hbox{$\pi\!:\cZ\!\lra\!\De$} of algebraic varieties that are locally of the~form
\BE{NCdegen_e} \pi\!:\big\{(\la,z_1,\ldots,z_k,p)\!\in\!\C^{k+2}\!\times\!\C^{n-k}\!:
z_1\!\ldots\!z_k\!=\!\la\big\}\lra\C,\quad
\pi\big(\la,z_1,\ldots,z_k,p\big)=\la,\EE
around the central fiber $\cZ_0\!\equiv\!\pi^{-1}(0)$.
The degenerations discussed above, i.e.~the standard one associated with
the symplectic sum construction of~\cite{Gf} and its extension indicated
in~\cite{SympDivConf,SympSumMulti}, correspond to $k\!=\!2$ in~\e_ref{NCdegen_e}.
The central fiber of~$\pi$ for $k\!\ge\!3$ in the algebro-geometric category
is a more general NC variety; see Figure~\ref{SympSum3_fig}.
Degeneration and splitting formulas for GW-invariants in this more general setting
require notions of GW-invariants for (smoothable) NC varieties and for smooth varieties
relative to NC divisors.
A degeneration formula in the projective category extending that of~\cite{JunLi2} 
has finally appeared in the setting of the \sf{logarithmic GW-theory} of~\cite{GS} in~\cite{ACGS};
the latter includes GW-invariants of smoothable NC varieties and of smooth varieties
relative to NC divisors. 
However, a splitting formula for the GW-invariants of NC varieties in the projective category 
remains to be established.\\

\noindent
The logarithmic GW-invariants of~\cite{GS} are special cases of the GW-invariants 
of exploded manifolds introduced in~\cite{BP-EX}.
Degeneration {\it and} splitting formulas for these invariants are 
studied in~\cite{BP}.
Based on the $k\!=\!2$ case established in~\cite{JunLi2},
one might expect that all curves in~$\cZ_0$ contributing to the GW-invariants of~$\cZ_{\la}$ 
either 
\begin{enumerate}[label=$\bu$,leftmargin=*]

\item have no irreducible components lying in the singular locus~$\cZ_0'$ of~$\cZ_0$
and meet at the smooth points of~$\cZ_0$ or at least 

\item have no irreducible components in~$\cZ_0'$.

\end{enumerate}
As demonstrated in~\cite{BP}, even the weaker alternative does not hold in general.
This makes any general splitting formula necessarily complicated;
its $k\!=\!3$ case is described in~\cite{BP3}.
A more geometric perspective on the GW-invariants of~\cite{BP-EX} appears in~\cite{Inc},
without analogues of the crucial degeneration and splitting formulas.
\\

\noindent
The GW-invariants of exploded manifolds of~\cite{BP-EX} and their interpretation in
some cases in~\cite{Inc} are essentially invariants of deformation equivalence
classes of almost K\"ahler structures on manifolds. 
While these classes are much larger than the deformation equivalence classes of
the algebro-geometric structures in \cite{GS,ACGS},
GW-invariants are fundamentally invariants of the still larger deformation equivalence
classes of symplectic structures.
Purely topological notions of \sf{NC symplectic divisors} and \sf{varieties}
are introduced in \cite{SympDivConf,SympDivConf2},
addressing a fundamental quandary of \cite[p343]{Gr} in the case of NC singularities.
Crucial to the introduction of these long desired notions is the new perspective proposed
in~\cite{SympDivConf}:\\

\begin{minipage}{6in}
{\it A symplectic variety/subvariety should be viewed as a deformation equivalence class
of objects with the same topology, not as a single object.}\\
\end{minipage}

\noindent
It is then shown in~\cite{SympDivConf,SympDivConf2} that the spaces of 
NC symplectic divisors and varieties are weakly homotopy equivalent to 
the spaces of almost K\"ahler structures, as needed for geometric applications.\\

\noindent
The equivalence between the topological and geometric notions of NC symplectic 
variety established in~\cite{SympDivConf,SympDivConf2} immediately implies
that any  invariants arising from~\cite{BP,Inc} in fact depend only on 
the deformation equivalence classes of symplectic structures.
These equivalences are also used in \cite{SympSumMulti,SympSumMulti2}
to establish a smoothability criterion for NC symplectic varieties.
Direct approaches to constructing GW-invariants of symplectic manifolds relative
to  NC symplectic divisors in the perspective of~\cite{SympDivConf,SympDivConf2} 
and to obtaining degeneration and splitting formulas for the degenerations appearing
in \cite{SympSumMulti,SympSumMulti2} are discussed in~\cite{MohNCrel} and~\cite{MohNCsum},
respectively.\\

\noindent
The decomposition and splitting formulas for GW-invariants in~\cite{BP}
involve \sf{exploded de~Rham cohomology} of~\cite{BPdeRm}, 
which makes these formulas very hard to apply.
The purpose of this elaborate modification of the ordinary de~Rham cohomology
is to correct the standard Kunneth decompositions of the diagonals of the strata
of the singular locus~$\cZ_0'$ of the central fiber~$\cZ_0$ for the presence of
lower-dimensional strata.
This removes certain \sf{degenerate contributions} to the Kunneth decompositions
of the diagonals of the strata of~$\cZ_0'$.
A local, completely topological approach to computing degenerate contributions 
in terms of the ordinary cohomology of the strata is presented in~\cite{g0pr}. 

\begin{ques}\label{NCdecomp_prob}
Is there a reasonably usable formula for general NC degenerations \hbox{$\pi\!:\cZ\!\lra\!\De$}
of symplectic manifolds which splits the GW-invariants of a smooth fiber~$\cZ_{\la}$ into
the GW-invariants of the strata of the central fiber~$\cZ_0$ that involves 
only the ordinary cohomology of the strata?
\end{ques}

\noindent
The introduction of symplectic topology notions of NC divisors, varieties, and degenerations
in \cite{SympDivConf,SympSumMulti,SympDivConf2,SympSumMulti2} has made it feasible
to study Question~\ref{NCdecomp_prob} entirely in the symplectic topology category,
which is far more flexible than the algebraic geometry category of~\cite{GS,ACGS}
and the almost K\"ahler category of~\cite{BP-EX,Inc}.
A symplectic approach to this question should fit well with the topological approach
of~\cite{g0pr} to degenerate contributions.
A splitting formula for GW-invariants of~$\cZ_{\la}$ resulting from such an approach
should involve sums over finite trees with the edges labeled by integer weights
and the vertices labeled by paths in the dual intersection complex of~$\cZ_0$
with additional de~Rham cohomology data;
these paths would correspond to the \sf{tropical curves} of~\cite{BP}.
While such a formula would still be more complicated than in the standard case of~\cite{JunLi2},
it should be more readily applicable than the presently available splitting formula
of~\cite{BP} that involves exploded de Rham cohomology.\\

\noindent
Degeneration and splitting formulas for real GW-invariants under real degenerations
of real symplectic manifolds have been obtained only in a small number of special cases.
A fundamental difficulty for obtaining such formulas is that the standard notions
of relative invariants of the complex GW-theory do not have direct analogues in
the real GW-theory in most settings.
Real GW-invariants of a real symplectic manifold $(X,\om,\phi)$ with simple contacts
with a real symplectic divisor $V\!\subset\!X$ can be readily defined 
whenever the real GW-invariants of $(X,\om,\phi)$ are defined and 
$V$ is disjoint from the fixed locus~$X^{\phi}$ of~$\phi$.
This observation lies behind the splitting formula and related vanishing result
for some genus~0 real GW-invariants under special real degenerations
of real symplectic manifolds obtained in~\cite{Teh1}.\\

\noindent
The reduction of the complex relative GW-invariants of $(X,V)$ to
the complex GW-invariants of~$X$ and~$V$ in~\cite{MauP} suggests the possibility
of expressing the real GW-invariants of a real symplectic sum $X_1\!\#_{X_{12}}\!X_2$
in terms of the real GW-invariants of $X_1,X_2,X_{12}$, whenever these are defined.
If $X_1$ and~$X_2$ are of real dimension~4, then $X$ is a real surface 
and the real GW-invariants of $X_1\!\#_{X_{12}}\!X_2$ should reduce to
the real GW-invariants of~$X_1$ and~$X_2$.
By \cite[Theorem~7]{BP15} and \cite[Theorem~1.1]{Brug16}, 
this is indeed the case for the genus~0 real GW-invariants 
if $X_{12}\!\approx\!\P^1$ is a real symplectic submanifold of 
self-intersection~2 in \hbox{$X_2\!=\!\P^1\!\times\!\P^1$}
and in some other settings with $X_{12}\!\approx\!\P^1$.
Genus~0 real GW-invariants have been defined for many real symplectic sixfolds and
for all real symplectic fourfolds. 
This leads to the following question.

\begin{ques}\label{RealSum_prob}
Is it possible to express the genus~0 real GW-invariants of 
a real symplectic sum $X_1\!\#_{X_{12}}\!X_2$ of 
real symplectic sixfolds $(X_i,\om_i,\phi_i)$ along 
a common real symplectic divisor~$X_{12}$ in terms of the genus~0 real GW-invariants 
$X_1,X_2,X_{12}$, whenever the genus~0 real GW-invariants  of the sixfolds are defined?
\end{ques}

\section{Geometric applications}
\label{GeomAppl_sec}

\noindent
Pseudoholomorphic curves were originally introduced in~\cite{Gromov}
with the aim of applications in symplectic topology.
These applications have included 
the Symplectic Non-Squeezing Theorem~\cite{Gromov},
classification of symplectic 4-manifolds~\cite{McDuff90,LaMc},
distinguishing diffeomorphic symplectic manifolds~\cite{Ruan94},
symplectic isotopy problem \cite{Ti98,SiTi},
and applications in birational algebraic geometry~\cite{Kollar,Zhiyu}.
However, many deep related problems remain open.\\

\noindent
Rational curves, i.e.~images of $J$-holomorphic maps from chains of spheres, 
play a particularly important role in algebraic geometry. 
A smooth algebraic manifold~$X$ is called \sf{uniruled} 
(resp.~\sf{rationally connected} or \sf{RC}) 
if there is a rational curve through every point 
(resp.~every pair of points) in~$X$.
According to~\cite{Kollar},
a uniruled algebraic variety admits a nonzero genus~0 GW-invariant
with a point insertion (i.e.~a count of stable maps in a fixed homology class
which pass through a point and some other constraints).
This implies that the uniruled property is invariant under symplectic deformations.
The RC property is known to be invariant under integrable deformations of
the complex structure \cite{Kollar}.
It is a long-standing conjecture of J.~Koll\'ar that the RC property is invariant under 
symplectic deformations as~well.
It is unknown if every RC algebraic manifold admits 
a nonzero genus~0 GW-invariant with two point insertions; 
this would immediately imply Koll\'ar's conjecture.
The dimension~3 case of this conjecture is established in~\cite{Zhiyu}
by combining 
the special cases treated in~\cite{Voisin} with the minimal model program.\\ 

\noindent
As GW-invariants are symplectic invariants, 
it is natural to consider the parallel situation in symplectic topology.
Given the flexibility of the symplectic category,
this may also provide a different approach to Koll\'ar's conjecture.
A symplectic manifold $(X,\om)$ is called \sf{uniruled} (resp.~{\sf{RC}})
if for some $\om$-compatible almost complex structure~$J$
 there is a genus~0 connected rational $J$-holomorphic curve through every point 
(resp.~every pair of points) in~$X$.
This leads to the following two pairs of questions.

\begin{ques}\label{uniruled_prob}
Let $J$ be {\it any} almost complex structure on
a uniruled (resp.~RC) compact symplectic manifold~$(X,\om)$.
Is there a connected rational $J$-holomorphic curve 
through every point (resp.~every pair of points) in~$X$?
\end{ques}

\begin{ques}\label{RCGW_prob}
Does every uniruled (resp.~RC) compact symplectic manifold~$(X,\om)$
admit a nonzero genus~0 GW-invariant with a point insertion (resp.~two point insertions)?
\end{ques}

\noindent
The affirmative answer to each case of Question~\ref{RCGW_prob} would immediately imply
the affirmative answer to the corresponding case of Question~\ref{uniruled_prob}.
The uniruled case of Question~\ref{RCGW_prob} is known only under the rigidity assumptions 
that $X$ is either K\"ahler \cite{Kollar} or admits a Hamiltonian $S^1$-action~\cite{McDuff}.
It is not difficult to construct $J$-holomorphic curves in a symplectic manifold that disappear
as the almost complex structure~$J$ deforms. 
On the other hand, regular $J$-holomorphic curves do not disappear under small deformations
of~$J$, while $J$-holomorphic curves contributing to nonzero GW-invariants 
survive all deformations of~$J$.
Thus, the above four questions concern the fundamental issue of 
the extent of flexibility in the symplectic category with implications 
to birational algebraic geometry.\\

\noindent
If $u\!:\P^1\!\lra\!X$ is a $J$-holomorphic map into a Kahler manifold
and for {\it some} $z\!\in\!\P^1$ the evaluation~map
\BE{evmap_e} H^0(\P^1;u^*TX)\lra T_{u(z)}X, \qquad \xi\lra\xi(z),\EE
is onto, then $H^1(\P^1;u^*TX)\!=\!0$, i.e.~$u$ is \textsf{regular}.
This statement is key to the arguments of~\cite{Kollar} in the algebraic setting.
It in particular implies that if the rational $J$-holomorphic curves cover
a nonempty open subset of a connected K\"ahler manifold,
then they cover all of~$X$.
As shown in~\cite{McLean16}, the last implication can fail in the almost K\"ahler category.
The first implication need not hold either, even if the evaluation homomorphism~\e_ref{evmap_e} 
is surjective for {\it every} $z\!\in\!\P^1$.
However,  the main results of~\cite{Kollar} may still extend to 
the almost K\"ahler category.
In particular, for the interplay between openness and closedness 
of various properties of complex structures exhibited in the proof of deformation 
invariance of the RC property for integrable complex structures in~\cite{Kollar}
to extend to a non-integrable complex structure,
the vanishing of the obstruction space needs to hold only generically
in a family of $J$-holomorphic maps covering~$X$.
This leads to a potentially even more fundamental problem in this spirit.

\begin{ques}\label{genreg_prob}
Let $\{u_{\al}\!:\P^1\!\lra\!X\}$ be a family of $J$-holomorphic curves
on a compact symplectic manifold~$(X,\om)$ that covers~$X$.
Is a generic member of this family a regular map?
\end{ques}

\noindent
There are still many open questions concerning the geography 
and topology of symplectic manifolds 
The multifold smoothing constructions of \cite{SympSumMulti,SympSumMulti2}
may shed light on some of these questions.
Just as the  ($2$-fold) symplectic sum construction of~\cite{Gf},
the multifold constructions could be used to build vast classes of non-K\"ahler symplectic 
manifolds with various topological properties.
They might also be useful for studying properties of symplectic manifolds 
of algebro-geometric flavor, in the spirit of the perspective on symplectic topology
initiated in~\cite{Gromov}.

\begin{ques}[{\cite[Question 14]{SympNCSumm}}]\label{RCsymp_ques}
Is every compact almost K\"ahler manifold with a rational 
$J$-holomorphic curve of a fixed homology class through every pair of points simply connected?
\end{ques}

\noindent
By \cite[Theorem~3.5]{Campana}, a compact RC K\"ahler manifold is simply connected.
As noted by J.~Starr, the fundamental group of a compact almost K\"ahler manifold $(X,\om,J)$
as in Question~\ref{RCsymp_ques} is finite.
The multifold sum/smoothing constructions of \cite{SympSumMulti,SympSumMulti2}
can be used to obtain symplectic manifolds that are not simply connected 
from simply connected ones and thus may be useful in answering Question~\ref{RCsymp_ques} 
negatively. 
The constructions of~\cite{SympSumMulti,SympSumMulti2} may also be useful in studying this 
question under the stronger assumption of
the existence of a nonzero GW-invariant of~$(X,\om)$ with two point insertions.\\

\noindent
As in the complex case, it is natural to expect that a real symplectic manifold $(X,\om,\phi)$
which has well-defined genus~0 real GW-invariants and is covered by real rational curves 
admits a nonzero genus~0 real GW-invariant with a real point insertion.
However, the reasoning neither in~\cite{Kollar},
which relies on the positivity of intersections in complex geometry,
nor in~\cite{McDuff}, which makes use of quantum cohomology,
is readily adaptable to the real setting.
Thus, there is not apparent approach to this problem at the present.\\

\noindent
Another important question in real algebraic geometry is the existence of real rational
curves on real even-degree complete intersections $X\!\subset\!\P^n$;
this would be implied by the existence of a well-defined nonzero genus~0 real GW-invariant
of~$X$.
However, the real analogue of the Quantum Lefschetz Hyperplane Principle~\e_ref{g0rel_e}
suggests that all such invariants should vanish.
On the other hand, one may hope for some real analogue of the reduced/family GW-invariants 
of \cite{BryanLeong,Junho}, which effectively remove a trivial line bundle from 
the obstruction cone for deformations of $J$-holomorphic maps to~$X$.
The resulting reduced/family real invariants could well be nonzero.\\

\noindent
{\it Department of Mathematics, Stony Brook University, Stony Brook, NY 11794\\
azinger@math.stonybrook.edu}


\begin{thebibliography}{99}
 

\bibitem{ACGS} D.~Abramovich, Q.~Chen, M.~Gross, and B.~Siebert, 
{\it Decomposition of degenerate Gromov-Witten invariants},
math/1709.09864

\bibitem{BCOV} M.~Bershadsky, S.~Cecotti, H.~Ooguri, and C.~Vafa,
{\it Holomorphic anomalies in topological field theories}, 
Nucl.~Phys.~B405 (1993), 279--304
 
\bibitem{Brug16} E.~Brugall\'e,
{\it Surgery of real symplectic fourfolds and Welschinger invariants},
math/1601.05708

\bibitem{BP15} E.~Brugall\'e and N.~Puignau, 
{\it On Welschinger invariants of symplectic 4-manifolds}, 
Comment.~Math.~Helv.~90 (2015), no.~4, 905–-938

\bibitem{BryanLeong} J.~Bryan and N.~Leung, 
{\it The enumerative geometry of K3 surfaces and modular forms},
J.~AMS 13 (2000), no.~2, 371–-410

\bibitem{BryanPand} J.~Bryan and R.~Pandharipande, 
{\it The local Gromov-Witten theory of curves}, J.~AMS 21 (2008), no.~1, 101-–136

\bibitem{CaDGP} P.~Candelas, X.~de la Ossa, P.~Green, and L.~Parkes,
{\it A Pair of Calabi-Yau manifolds as an exactly soluble superconformal theory},
Nuclear Phys.~B359 (1991), 21--74

\bibitem{Campana} F.~Campana, 
{\it On twistor spaces of the class $\cC$},
J.~Diff.~Geom.~33 (1991), no.~2, 541-–549

\bibitem{CH} L.~Caporaso and J.~Harris,
{\it Counting plane curves in any genus}, Invent.~Math.~131 (1998), no.~2, 345--392

\bibitem{DoanWalpulski} A.~Doan and T.~Walpulski,
{\it Equivariant Brill-Noether theory for elliptic operators and 
applications to symplectic geometry}, preprint

\bibitem{Eft} E.~Eftekhary, 
{\it On finiteness and rigidity of $J$-holomorphic curves in symplectic three-folds}.
Adv.~Math.~289 (2016), 1082–-1105

\bibitem{ES} G.~Ellingsrud and S.~Str\"omme,
{\it Bott's formula and enumerative geometry},
J.~AMS 9 (1996), no.~1, 175--193

\bibitem{Teh1} M.~Farajzadeh Tehrani,
{\it Open Gromov-Witten theory on symplectic manifolds and symplectic cutting},
Adv.~Math.~232 (2013), no.~1, 238–-270

\bibitem{MohNCrel} M.~Farajzadeh Tehrani,
{\it Pseudoholomorphic maps relative to normal crossings symplectic divisors: compactification},
math/1710.00224

\bibitem{MohNCsum} M.~Farajzadeh Tehrani,
{\it Towards a degeneration formula for the Gromov-Witten invariants of symplectic manifolds},
math/1710.00599

\bibitem{SympNCSumm} M.~Farajzadeh Tehrani, M.~McLean, and A.~Zinger, 
{\it Singularities and semistable degenerations for symplectic topology},
C.~R.~Math. Acad.~Sci.~Paris 356 (2018), no.~4, 420–-432

\bibitem{SympDivConf} M.~Farajzadeh Tehrani, M.~McLean, and A.~Zinger,
{\it Normal crossings singularities for symplectic topology},
math/1410.0609v3

\bibitem{SympSumMulti} M.~Farajzadeh Tehrani, M.~McLean, and A.~Zinger,
{\it The smoothability of normal crossings symplectic varieties}, math/1707.01464

\bibitem{SympDivConf2} M.~Farajzadeh-Tehrani, M.~McLean, and A.~Zinger,
{\it Normal crossings singularities for symplectic topology, II}, in preparation

\bibitem{SympSumMulti2} M.~Farajzadeh-Tehrani, M.~McLean, and A.~Zinger,
{\it The smoothability of normal crossings symplectic varieties, II}, in preparation

\bibitem{SympSum} M.~Farajzadeh Tehrani and A.~Zinger,
{\it On symplectic sum formulas in Gromov-Witten theory}, 
math/1404.1898

\bibitem{GWrelIP} M.~Farajzadeh Tehrani and A.~Zinger,
{\it On the rim tori refinement of relative Gromov-Witten invariants},
math/1412.8204

\bibitem{GWsumIP} M.~Farajzadeh Tehrani and A.~Zinger,
{\it On the refined symplectic sum formula for Gromov-Witten invariants},
math/1412.8205

\bibitem{FO} K.~Fukaya and K.~Ono,
{\it Arnold Conjecture and Gromov-Witten invariant},
Topology 38 (1999), no.~5, 933--1048

\bibitem{FP} W.~Fulton and R.~Pandharipande, 
{\it Notes on stable maps and quantum cohomology},
Proc.~Sympos.~Pure Math.~62, Part~2, 45–-96, AMS, 1997

\bibitem{Ge2} P.~Georgieva,
{\it Open Gromov-Witten invariants in the presence of an anti-symplectic involution},
Adv.~Math.~301 (2016), 116–-160

\bibitem{RealGV} P.~Georgieva and E.~Ionel, seminar presentations~2017,~2018

\bibitem{RealEnum} P.~Georgieva and A.~Zinger,
{\it Enumeration of real curves in $\C\P^{2n-1}$ and
a WDVV relation for real Gromov-Witten invariants}, 
Duke Math.~J.~166 (2017), no.~17, 3291–-3347

\bibitem{RealEnumApp} P.~Georgieva and A.~Zinger,
{\it A recursion for counts of real curves in $\C\P^{2n-1}$: another proof},
 Internat.~J.~Math. 29 (2018), no.~4, 1850027

\bibitem{RealGWsI} P.~Georgieva and A.~Zinger,
{\it Real Gromov-Witten theory in all genera and real enumerative geometry: construction},
math/1504.06617v5

\bibitem{RealGWsII} P.~Georgieva and A.~Zinger,
{\it Real Gromov-Witten theory in all genera and real enumerative geometry: properties},
math/1507.06633v4

\bibitem{RealGWsIII} P.~Georgieva and A.~Zinger,
{\it Real Gromov-Witten theory in all genera and real enumerative geometry: computation},
math/1510.07568

\bibitem{Getz} E.~Getzler, 
{\it Intersection theory on $\ov\cM_{1,4}$ and elliptic Gromov-Witten invariants}, 
J.~AMS 10 (1997), no.~4, 973--998

\bibitem{g0ms} A.~Givental, 
{\it The mirror formula for quintic threefolds},
{\it California Symplectic Geometry Seminar}, 
AMS Transl.~Ser.~196, 49--62, AMS, 1999

\bibitem{Giv} A.~Givental, 
{\it Gromov-Witten invariants and quantization of quadratic Hamiltonians},
Mosc.~Math.~J.~1 (2001), no.~4, 551–-568

\bibitem{GH} P.~Griffiths and J.~Harris, 
{\it Principles of Algebraic Geometry}, Wiley Classics Library,~1994

\bibitem{Gf} R.~Gompf, {\it A new construction of symplectic manifolds}, 
Ann.~of Math.~142 (1995), no.~3, 527--595

\bibitem{GV} R.~Gopakumar and C.~Vafa, 
{\it M-theory and topological strings I,II}, 
hep-th/9809187,9812127

\bibitem{Gromov} M.~Gromov,
{\it Pseudoholomorphic curves in symplectic manifolds},
Invent.~Math.~82 (1985), no.~2, 307--347

\bibitem{Gr} M.~Gromov, {\it Partial Differential Relations}, Springer-Verlag, 1986

\bibitem{GS0} M.~Gross and B.~Siebert,
{\it  Affine manifolds, log structures, and mirror symmetry}, 
Turkish J.~Math.~27 (2003), no.~1, 33–-60

\bibitem{GS} M.~Gross and B.~Siebert, 
{\it Logarithmic Gromov-Witten invariants}, J.~AMS 26 (2013), no.~2, 451--510

\bibitem{Inc} E.~Ionel,
{\it GW-invariants relative normal crossings divisors},
Adv.~Math.~281 (2015), 40-–141

\bibitem{IPrel} E.~Ionel and T.~Parker, 
{\it Relative Gromov-Witten invariants},
Ann.~of Math.~157 (2003), no.~1, 45--96

\bibitem{IPsum} E.~Ionel and T.~Parker, 
{\it The symplectic sum formula for Gromov-Witten invariants},
Ann.~of Math.~159 (2004), no.~3, 935--1025

\bibitem{IPgv} E.~Ionel and T.~Parker, 
{\it The Gopakumar-Vafa formula for symplectic manifolds}, 
Ann.~of Math.~187 (2018), no.~1, 1-–64

\bibitem{KlemmP} A.~Klemm and R.~Pandharipande,
{\it Enumerative geometry of Calabi-Yau 4-folds},
Comm.~Math.~Phys.~261 (2006), no.~2, 451--516

\bibitem{Kollar} J.~Koll\'ar, 
{\it Rational Curves on Algebraic Varieties}, 
Modern Surveys in Mathematics, Springer-Verlag, 1996

\bibitem{Kont} M.~Kontsevich, 
{\it Enumeration of rational curves via torus actions}, 
The Moduli Space of Curves,  Progr.~Math.~129 (1995), 335–-368

\bibitem{KM} M.~Kontsevich and Y.~Manin, 
{\it Gromov-Witten classes, quantum cohomology, and enumerative geometry}, 
Comm.~Math.~Phys.~164 (1994), no.~3, 525--562

\bibitem{LaMc} F.~Lalonde and D.~McDuff, 
{\it $J$-curves and the classification of rational and ruled symplectic 4-manifolds},
Contact and Symplectic Geometry, 3–-42, Publ.~Newton Inst.~8, 1996

\bibitem{LM} B.~Lawson and M.-L.~Michelsohn, 
{\it Spin Geometry}, Princeton Mathematical Series~38, 1989

\bibitem{Junho} J.~Lee, 
{\it Family Gromov-Witten invariants for K\"ahler surfaces}, 
Duke Math.~J.~123 (2004), no.~1, 209-–233

\bibitem{JunLi1} J.~Li, 
{\it Stable morphisms to singular schemes and relative stable morphisms}, 
J.~Diff.~Geom.~57 (2001), no.~3, 509--578

\bibitem{JunLi2} J.~Li, {\it A degeneration formula of GW-invariants}, 
J.~Differential Geom.~60 (2002), no.~3, 199–-293

\bibitem{LT0} J.~Li and G.~Tian, 
{\it Virtual moduli cycles and Gromov-Witten invariants of algebraic varieties}, 
J.~AMS 11 (1998), no.~1, 119-–174

\bibitem{LT} J.~Li and G.~Tian, 
{\it Virtual moduli cycles and Gromov-Witten invariants of general symplectic manifolds}, 
Topics in Symplectic \hbox{$4$-Manifolds},
47-83, First Int.~Press Lect.~Ser., I, Internat.~Press, 1998

\bibitem{g1gw} J.~Li and A.~Zinger, 
{\it On the genus-one Gromov-Witten invariants of complete intersections}, 
J.~Diff.~Geom.~82 (2009), no.~3, 641--690

\bibitem{LLY} B.~Lian, K.~Liu, and S.-T.~Yau, 
{\it Mirror Principle I}, Asian J.~of Math.~1, No.~4 (1997), 729--763

\bibitem{XiaoboTian} X.~Liu and G.~Tian,
{\it Virasoro constraints for quantum cohomology}, 
J.~Diff.~Geom.~50 (1998), no.~3, 537-–590

\bibitem{MauP} D.~Maulik and R.~Pandharipande, 
{\it A topological view of Gromov-Witten theory}, 
Topology~45 (2006), no.~5, 887-–918

\bibitem{MauP2} D.~Maulik and R.~Pandharipande, 
{\it New calculations in Gromov-Witten theory},
Pure Appl.~Math.~Q.~4 (2008), no.~2, 469–-500

\bibitem{McDuff90} D.~McDuff, 
{\it The structure of rational and ruled symplectic four-manifolds}, 
J.~AMS 3 (1990), 679-–712

\bibitem{McDuff} D.~McDuff,
{\it Hamiltonian $S^1$-manifolds are uniruled}, 
Duke Math.~J.~146 (2009), no.~3, 449--507

\bibitem{McLean16} M.~McLean,
{\it Local uniruledness does Not imply global uniruledness
in symplectic topology}, available from the author's website

\bibitem{NakYos} H.~Nakajima and K.~Yoshioka, 
{\it Instanton counting on blowup,~I: 4-dimensional pure gauge theory}, 
Invent.~Math.~162 (2005), no.~2, 313-–355

\bibitem{g2comp} J.~Niu,
{\it A sharp compactness theorem for genus-two pseudo-holomorphic maps},
PhD thesis, Stony Brook,~2016

\bibitem{NZ} J.~Niu and A.~Zinger,
{\it Lower bounds for the enumerative geometry of positive-genus real curves},
math/1511.02206

\bibitem{OP1} A.~Okounkov and R.~Pandharipande,
{\it Virasoro constraints for target curves},  
Invent.~Math.~163  (2006),  no.~1, 47--108, AMS, 2009

\bibitem{OP0} A.~Okounkov and R.~Pandharipande,
{\it Gromov-Witten theory, Hurwitz numbers, and matrix models},
Proc.~Sympos.~Pure Math.~80, Part 1, 325–-414

\bibitem{PandQdiv}  R.~Pandharipande, 
{\it Intersections of $\Q$-divisors on Kontsevich's moduli space 
$\ov{M}_{0,n}(\P^r,d)$ and enumerative geometry}, 
Trans.~AMS 351 (1999), no.~4, 1481--1505

\bibitem{Pand99} R.~Pandharipande, 
{\it A geometric construction of Getzler's elliptic relation}, 
Math.~Ann. 313 (1999), no.~4, 715--729

\bibitem{Pand3ques} R.~Pandharipande, 
{\it Three questions in Gromov-Witten theory}, 
Proceedings of ICM, Vol.~II, 503–-512, 2002

\bibitem{cy5fold} R.~Pandharipande and A.~Zinger, 
{\it Enumerative Geometry of Calabi-Yau 5-folds},
New Developments in Algebraic Geometry, Integrable Systems, and Mirror Symmetry, 
239--288,  Adv.~Stud.~Pure Math.~59, 2010 

\bibitem{Pardon} J.~Pardon, 
{\it An algebraic approach to virtual fundamental cycles on moduli spaces 
of pseudo-holomorphic curves}, Geom.~Topol.~20 (2016), no.~2, 779-–1034

\bibitem{BP-EX}  B.~Parker, {\it Exploded fibrations}, 
Proceedings of Gokova Geometry-Topology Conference 2006 (2007), 52--90

\bibitem{BP} B.~Parker, {\it Gromov-Witten invariants of exploded manifolds},
math/1102.0158

\bibitem{BP3} B.~Parker,
{\it Gluing formula for Gromov-Witten invariants in a triple product},
math/1511.00779

\bibitem{BPdeRm} B.~Parker, {\it de Rham theory of exploded manifolds},
Geom.~Topol.~22 (2018), no.~1, 1-–54

\bibitem{Pe14} D.~Petersen, 
{\it The structure of the tautological ring in genus one},
Duke Math.~J.~163 (2014), no.~4, 777-–793 

\bibitem{Pe16} D.~Petersen, 
{\it Tautological rings of spaces of pointed genus two curves of compact type}, 
Compos.~Math.~152 (2016), no.~7, 1398-–1420 

\bibitem{PeTo} D.~Petersen and O.~Tommasi,
{\it The Gorenstein conjecture fails for the tautological ring of $\ov\cM_{2,n}$},
Invent.~Math.~196 (2014), no.~1, 139-–161

\bibitem{Po} A.~Popa, 
{\it The genus one Gromov-Witten invariants of Calabi-Yau complete intersections},
Trans.~AMS 365 (2013), no.~3, 1149--1181

\bibitem{RiTi} T.~Rivi\'ere and G.~Tian, 
{\it The singular set of 1-1 integral currents},
Ann.~of Math.~169 (2009), no.~3, 741-–794

\bibitem{Ruan94} Y.~Ruan, 
{\it Symplectic topology on algebraic 3-folds}, 
J.~Diff.~Geom.~39 (1994), no.~1, 215--227

\bibitem{RT} Y.~Ruan and G.~Tian, 
{\it A mathematical theory of quantum cohomology}, 
J.~Diff.~Geom.~42 (1995), no.~2, 259--367

\bibitem{RT2} Y.~Ruan and G.~Tian, 
{\it Higher genus symplectic invariants and sigma models coupled with gravity}, 
Invent.~Math.~130 (1997), no.~3, 455--516

\bibitem{SiTi} B.~Siebert and G.~Tian, 
{\it On the holomorphicity of genus two Lefschetz fibrations}, 
Ann.~of Math.~161 (2005), no.~2, 959–-1020

\bibitem{Sol2} J.~Solomon, 
{\it A differential equation for the open Gromov-Witten potential},
 pre-print 2007

\bibitem{Tel} C.~Teleman, 
{\it The structure of 2D semi-simple field theories},
Invent.~Math.~188 (2012), no.~3, 525-–588, 2012

\bibitem{Ti} G.~Tian, 
{\it The quantum cohomology and its associativity}, 
Current Developments in Mathematics~1995, 361--401, Inter.~Press.

\bibitem{Ti98} G.~Tian,
{\it Symplectic isotopy in four dimension}, 
First International Congress of Chinese Mathematicians, 143-–147, 
AMS/IP Stud.~Adv.~Math.~20, 2001

\bibitem{Zhiyu} Z.~Tian, {\it Symplectic geometry of rationally connected
threefolds},  Duke Math.~J.~161 (2012), no.~5, 803--843

\bibitem{g1desing0} R.~Vakil and A.~Zinger, 
{\it A natural smooth compactification of the space of elliptic curves in
projective space}, ERA AMS 13 (2007), 53--59

\bibitem{g1desing} R.~Vakil and A.~Zinger, 
{\it A desingularization of the main component of
the moduli space of genus-one stable maps into $\P^n$}, 
Geom.~Topol.~12 (2008), no.~1, 1--95

\bibitem{Voisin} C.~Voisin, {\it Rationally connected 3-folds and symplectic 
geometry}, Ast\'erisque 322 (2008), 1--21

\bibitem{Wal} J.~Walcher,   
{\it Evidence for tadpole cancellation in the topological string}, 
Comm.~Number Theory Phys.~3 (2009), no.~1, 111--172

\bibitem{Wel4}  J.-Y.~Welschinger,  
{\it Invariants of real symplectic 4-manifolds and lower bounds in real enumerative geometry},
Invent.~Math.~162 (2005), no.~1, 195--234

\bibitem{Wendl} C.~Wendl,
{\it Transversality and super-rigidity for multiply covered holomorphic curves},
math/1609.09867

\bibitem{g0pr} A.~Zinger, 
{\it Counting rational curves of arbitrary shape in projective spaces}, 
Geom.~Topol.~9 (2005), 571-–697

\bibitem{g1cone} A.~Zinger, 
{\it On the structure of certain natural cones over moduli spaces of 
genus-one holomorphic maps}, Adv.~Math.~214 (2007), no.~2, 878--933

\bibitem{pseudo} A.~Zinger, {\it Pseudocycles and integral homology},
Trans.~AMS 360 (2008), no.~5, 2741--2765

\bibitem{g1diff} A.~Zinger, 
{\it Standard vs.~reduced genus-one Gromov-Witten invariants},
Geom.~Topol.~12 (2008), no.~2, 1203--1241

\bibitem{g1comp} A.~Zinger,  
{\it A sharp compactness theorem for genus-one pseudo-holomorphic maps},  
Geom.~Topol.~13 (2009), no.~5, 2427--2522

\bibitem{g1comp2} A.~Zinger, {\it Reduced genus-one Gromov-Witten invariants},
J.~Diff.~Geom.~83 (2009), no.~2, 407--460

\bibitem{bcov1} A.~Zinger,
{\it The reduced genus 1 Gromov-Witten invariants of Calabi-Yau hypersurfaces},
J.~AMS~22 (2009), no.~3, 691–-737

\bibitem{FanoGV} A.~Zinger, 
{\it A comparison theorem for Gromov-Witten invariants in the symplectic category},
Adv.~Math.~228 (2011), no.~1, 535--574

\end{thebibliography}
\end{document}